\documentclass[a4paper]{article}

\usepackage{amsthm,amsmath,amssymb}
\usepackage{mathrsfs}
\usepackage{graphicx}
\usepackage{float,verbatim}
\usepackage{threeparttable,booktabs}
\usepackage{underscore}
\usepackage{algorithm}
\usepackage{subfigure}
\usepackage[noend]{algpseudocode}
\usepackage[nohead,margin=1.0in]{geometry}
\usepackage{color}
\usepackage{hyperref}
\hypersetup{hypertex=true,
	colorlinks=true,
	linkcolor=blue,
	anchorcolor=blue,
	citecolor=blue}
\graphicspath{{./pic/}}

\newtheorem{proposition}{Proposition}[section]
\newtheorem{corollary}{Corollary}[section]
\newtheorem{theorem}{Theorem}[section]
\newtheorem{remark}{Remark}[]
\newtheorem{lemma}{Lemma}[section]
\newtheorem{example}{Example}[section]

\numberwithin{equation}{section}

\title{Solving Elliptic Problems with Singular Sources using Singularity Splitting Deep Ritz Method\thanks{The work of B. Jin is supported by UK EPSRC grant EP/T000864/1 and EP/V026259/1, and a start-up fund from The Chinese University of Hong Kong. The work of  Z. Zhou is supported by Hong Kong Research Grants Council (15303122) and an internal grant of Hong Kong Polytechnic University (Project ID: P0031041, Work Programme: ZZKS).}}
\author{Tianhao Hu\thanks{School of Mathematics, Jilin University, Changchun 130012, P.R. China (\texttt{huth1019@mails.jlu.edu.cn})}
\and Bangti Jin\thanks{Department of Mathematics, The Chinese University of Hong Kong, Shatin, New Territories, Hong Kong, P.R. China (\texttt{bangti.jin@gmail.com, b.jin@
cuhk.edu.hk}).} \and Zhi Zhou\thanks{Department of Applied Mathematics,
The Hong Kong Polytechnic University, Kowloon, Hong Kong, P.R. China (\texttt{zhizhou@polyu.edu.hk})}}
\date{}

\begin{document}
	\maketitle
	
\begin{abstract}
In this work, we develop an efficient solver  based on neural networks for second-order elliptic  equations with variable coefficients and singular sources.
This class of problems covers general point sources, line sources and the combination of point-line sources, and has a broad range of practical applications. The proposed approach is based on decomposing the true solution into a singular part that is known analytically using the fundamental solution of the Laplace equation and a regular part that satisfies a suitable modified elliptic PDE with a smoother source, and then solving for the regular part using the deep Ritz method. A path-following strategy is suggested to select the penalty parameter for enforcing the Dirichlet boundary condition. Extensive numerical experiments in two- and multi-dimensional spaces with point sources, line sources or their combinations are presented to illustrate the efficiency of the proposed approach, and a comparative study with several existing approaches based on neural networks is also given, which shows clearly its competitiveness for the specific class of problems. In addition, we briefly discuss the error analysis of the approach.
\end{abstract}

\section{Introduction}

Partial differential equations (PDEs) represent a very important class of mathematical models that plays a vital role in physics, science and engineering. Traditional numerical methods for solving PDEs include finite difference method (FDM), finite element method (FEM), finite volume method (FVM) and boundary element method (BEM) etc. These methods have been maturely developed over the past few decades and efficient implementations and rigorous theoretical guarantees, e.g., error estimates, are also readily available. In recent years, neural networks (NNs) have emerged as a promising way to solve PDEs \cite{lagaris1998artificial}, motivated by their great successes in computer vision, speech recognition and natural language processing etc. A large number of neural PDE solvers have been developed, including physics informed neural networks (PINNs) \cite{RAISSI2019686}, deep Ritz method (DRM) \cite{yu2018deep}, deep Galerkin method (DGM) \cite{sirignano2018dgm}, weak adversarial network (WAN) \cite{Zang:2020} and deep least-squares method \cite{CaiChenLiu:2021}, to name just a few. Compared with more traditional methods, deep learning solvers have demonstrated very encouraging results in several direct and inverse problems \cite{Karniadakis:2021nature,EHanJentzen:2022,JinLiLu:2022ip,CuomoSchiano:2022}.

Note that in all these approaches, one uses NNs as ansatz functions to approximate the solutions to the PDEs. According to existing approximation theory of NNs \cite{DeVoreHanin:2021}, the convergence rates of these methods depend heavily on the regularity of the solution (as well as the stability of the mathematical formulation). A direct application of these methods might be ineffective or even fail spectacularly when dealing with challenging scenarios \cite{WangPerdikaris:2022jcp,Krishnapriyan:2021}, e.g., convection-dominated problems, transport problem, high-frequency wave propagation, problems with geometric singularities (cracks / corner singularity) and singular sources. Note that all these settings lead to either strong directional behavior or weak solution singularities or highly oscillatory behavior, which are challenging for physics-agnostic NNs to capture effectively.

\subsection{Mathematical formulation}\label{sec:formulation}

This work deals with one challenging class of boundary value problems, i.e., second-order elliptic problems of divergence form with singular sources. Let $\Omega$ be an open bounded domain in $\mathbb{R}^d$ ($d\geq 2$) with a smooth boundary $\partial\Omega$.
Without loss of generality, consider the following variable coefficient second-order elliptic equation in the divergence form with a singular source:
\begin{equation}\label{problem}
		\begin{aligned}
			-\nabla\cdot (\kappa\nabla u)=S+g,\quad \mbox{in }\Omega,
		\end{aligned}
	\end{equation}
where the scalar-valued diffusion coefficient $\kappa\in C^2(\overline\Omega)$ satisfies uniform ellipticity, i.e., $\kappa(\mathbf{x}) \geq m$, for some fixed constant $m>0$, and $S$ and $g$ denote a singular and smooth source, respectively. Here the singular source $S$ can be expressed in terms of the Dirac delta function, and covers both point and line sources. Specifically, let $\delta(\mathbf{x})$ denote the Dirac delta function concentrated at the origin, defined in the sense of distribution, i.e.,
$\langle\delta, v \rangle = v(\mathbf{0})$ for all $v\in C(\overline{\Omega})$. Likewise, with $\Lambda$ being
a line / line segment in the domain $\Omega$, let $\delta_\Lambda$ be a Dirac measure concentrated
on $\Lambda$, such that $\langle\delta_\Lambda, v \rangle = \int_\Lambda v \, \d s$ for all $v\in C(\overline{\Omega})$.
Then the singular source $S$ can take the form
\begin{equation*}
    S(\mathbf{x}) = \sum_{i=1}^{N_p} c_i\delta (\mathbf{x}-\mathbf{x}_i) + \sum_{j=1}^{N_L}f_j\delta_{\Lambda_j}(\mathbf{x}),
\end{equation*}
where $\{(c_i,\mathbf{x}_i): c_i\in\mathbb{R},\mathbf{x}_i\in\Omega\}_{i=1}^{N_p}$ denotes a collection of $N_p$ Dirac
point sources in the domain, concentrated at $\mathbf{x}_i$ with a strength $c_i$, and
$\{(f_j,\Lambda_j): f_j\in C(\overline{\Omega}), \,\Lambda_j\subset\Omega\}_{j=1}^{N_L}$denotes a collection
of $N_L$ line sources, supported on the line segment $\Lambda_j$ with the corresponding density $f_j$. The proposed approach
can be adapted to surface sources etc.
Singular sources supported on high-dimensional sets induce a weaker singularity in the solution than point sources.

This class of problems has a broad range of practical applications. For example, a point source can describe a pulse exciting the electric field in electromagnetic simulation \cite{evans10,sullivan2013electromagnetic}, or volume injection in an acoustic wave equation \cite{moseley2020solving}.
In practice, the source can concentrate not only on points but also on lines, curves or surfaces, and the latter arises naturally in biological modeling. For example, the works \cite{grinberg2011modeling,reichold2009vascular,moore2005one} discuss blood flow in the vascularized tissue of the brain, and \cite{cattaneo2014computational,d2008coupling} studied drug delivery through microcirculation and tissue perfusion. In these applications, the resulting mathematical models are PDEs with line source(s). Hence it is also of great interest to develop efficient numerical solvers for these problems.

Elliptic problems with point sources have been extensively
studied in the context of classical numerical PDE solvers, and there are several different approaches to handle the
singularity. In some situations, traditional solvers, mostly FEM, can be applied directly, but require careful interpretation and several (weighted) $L^2(\Omega)$ error estimates have been obtained (see, e.g., \cite{Babuska:1971,Apel:2011,2012FINITE,KoppWohlmuth:2014}).
Due to the limited solution regularity, there is no error estimate in the $ H^1 $ norm. Alternative, one may regularize the Dirac
function with a smooth function (e.g., Gaussians or Taylor expansion technique), or to approximate the solution in weighted Sobolev spaces (with a weight vanishing
at the singularity) so that the singularity of the solution plays a less prominent role  (see \cite{HOSSEINI2016423,
ENGQUIST200528,liu2012properties} and references therein).
Once the source singularity is regularized, the problem can then be solved using any stand-alone techniques, e.g., FDM or FEM. There are also several alternative approaches taking into account the problem structures, e.g., mesh grading \cite{Apel:2011}, generalized FEM \cite{Fries:2010} and singularity reconstruction \cite{cai2001finite}. In particular, these approaches build the analytical insights into the construction of the numerical procedures directly in order to achieve better computational efficiency. So far most studies are concerned with point sources. The case of line sources is far less studied. The few exceptions include D'Angelo and Quarteroni \cite{d2008coupling,2012FINITE} which gave the existence of a solution to the variational formulation (in weighted Sobolev spaces) and studied its finite element approximation, and Gjerde et al \cite{gjerde2019splitting,GjerdeNordbotten:2021sinum}, where a novel singularity splitting technique was suggested and investigated within the context of the Galerkin and mixed formulations. In particular, the singularity splitting technique \cite{gjerde2019splitting,GjerdeNordbotten:2021sinum} enjoys improved approximation properties: it yields optimal convergence rates for lowest-order elements and removes the pollution around the line source \cite{gjerde2019splitting}.

NN based solvers have not been extensively studied for this class of problems. Due to the low-regularity of the source term, the standard PINN or deep Ritz formulations cannot be applied directly. Indeed, the corresponding mathematical formulations are not always well-defined, and this is also clearly manifested in the numerical experiments. See Section \ref{ssec:existing} for detailed discussions on this aspect of the challenges. It is worth pointing out that due to the persistent optimization error, none of existing  neural PDE solvers can exhibit a steady and consistent convergence rate in practical computation, which is in stark contrast to classical PDE solvers, e.g., FEM and FDM. See the recent works \cite{AdcockDexter:2021, Cyr:2020} for detailed discussions on this important issue.

\subsection{The contributions of this work and organization}
In this work we propose a modification of the DRM \cite{yu2018deep}, called singularity splitting DRM (SSDRM), that can be applied effectively for the numerical solution of second-order elliptic PDEs with singular sources.
 Our approach is inspired by \cite{gjerde2019splitting,10.21186893PA}: we split the solution into a singular
and a regular part, and then approximate only the regular part using NNs. The singular part is explicitly expressed in terms of the fundamental solution $\Phi$ to the Laplace equation, and the regular part is smoother and satisfies a modified elliptic PDE that can be effectively solved using the DRM. Due to the improved
regularity of the regular part, NNs are effective for the numerical approximation, as indicated by the approximation
theory of NNs \cite{yarotsky2017error,I2020Approximation}.
The overall approach is as easy to implement as the original DRM, but significantly improves the accuracy of the approximation by carefully exploiting the known analytic structure of the underlying exact solution, which follows closely the current paradigm of building physical insights into machine learning techniques for problems in computational science and engineering.

The extensive numerical experiments in two- and multi-dimension in Section \ref{sec:numer} show clearly the efficiency and accuracy of the proposed SSDRM, especially extracting the solution singularity can greatly enhance the approximation accuracy. The comparative study shows that it significantly outperforms existing NN based approaches in terms of accuracy. In addition, we also briefly discuss the associated theoretical issues of the approach, which provide guidelines on the choice of the penalty parameter and the neural network architecture. More precisely, our contributions can be summarized as:
\begin{itemize}
    \item[(i)] to develop a novel singularity splitting based technique for problem \eqref{problem} to enable the use of DRM.
    \item[(ii)] to provide relevant theoretical underpinning of SSDRM, including both penalization and generalization errors.
\item[(iii)] to carry out extensive numerical examples  in two- and multi-dimension,
with one or multiple singularities to illustrate the flexibility and accuracy of SSDRM.
\item[(iv)] to conduct a comparative study with several state-of-the-art NN based PDE solvers.
\end{itemize}

The rest of the paper is organized as follows. In Section \ref{ssec:existing} we review existing NN based PDE solvers, and discuss challenges associated with their applications to problem \eqref{problem}. In Section \ref{sec:method} we propose the numerical method, based on the idea of singularity splitting. Then in Section \ref{sec:numer} we present extensive numerical experiments to illustrate the accuracy and efficiency of the approach. Finally we discuss relevant theoretical issues in Section \ref{sec:theory}, and conclude the paper with a short conclusion.
Throughout, we denote by $(\cdot,\cdot)$ and $(\cdot,\cdot)_{L^2(\partial\Omega)}$ the $L^2(\Omega)$ and $L^2(\partial\Omega)$ inner products, respectively. The notation $|\cdot|$ and $\cdot$ denote respectively the usual Euclidean norm and inner product on $\mathbb{R}^d$.

\section{Challenges for Neural Network Based Solvers}\label{ssec:existing}
The presence of a singular source in the model \eqref{problem} poses big challenges to its efficient numerical approximation. This is particularly severe for NN based approaches, which tend to approximate smooth solutions well due to their implicit structural bias \cite{Soudry:2018,WangPerdikaris:2021sisc,WangPerdikaris:2022jcp}, but not weakly singular solutions as arising in problem \eqref{problem}. Below we illustrate the inherent computational challenges for three existing NN solvers, i.e., PINN, DRM and WAN, for solving a model version of problem \eqref{problem} with a Dirichlet boundary condition
\begin{align}\label{diri}
  \left\{\begin{aligned}
  -\nabla\cdot(\kappa\nabla u) &= \delta(\mathbf{x-x}_0)+g,&& \mbox{in }\Omega,\\
    u&=h,&& \mbox{on }\partial\Omega,
  \end{aligned}\right.
\end{align}
(i.e., with a point source $\delta(\mathbf{x}-\mathbf{x}_0)$, $\mathbf{x}_0\in\Omega$) in their continuous formulations. In practice, all these solvers approximate the PDE solution $u$ directly using NNs, which are then discretized properly.

Physics Informed Neural Network (PINN) is based on PDE residual minimization. When applying PINN to solve problem \eqref{diri} directly, the loss function $L(u)$ is commonly taken to be
\begin{equation}
\min_{u\in H^1(\Omega)} L(u)=L_r(u)+\lambda L_{b}(u),
\end{equation}
where $\lambda>0$ is a penalty parameter to approximately enforce the Dirichlet boundary condition in \eqref{diri}, and the physical loss $L_r(u)$ and boundary loss $L_b(u)$ are given respectively by
\begin{align*}
    L_r(u)=\|-\nabla\cdot (\kappa\nabla u)-\delta(\mathbf{x}-\mathbf{x}_0)-g\|_{L^2(\Omega)}^2\quad\mbox{and}\quad
   L_{b}(u)=\|u-h\|_{L^2(\partial\Omega)}^2.
\end{align*}
Note that the point source $\delta(\mathbf{x}-\mathbf{x}_0)$ does not belong to the space $L^2(\Omega)$ but only the space of Radon measures, and thus the standard PINN loss $L_r(u)$ with the $L^2(\Omega)$ norm of the PDE residual is not well defined. Hence, one cannot apply PINN directly to problem \eqref{diri}, and this is also confirmed by numerical experiments. Indeed, when applying PINN loss to solve \eqref{diri} naively, there are pronounced errors near the singularity \cite{rasht2022physics, moseley2020solving}. This calls for proper modification of the original PINN in order to obtain good approximations, e.g., self-adaptive PINN described below. There are alternative formulations of PINN, in the spirit of the weak formulation for PDEs, e.g., variational PINN \cite{KharazmiZhangKarniadakis:2021}, that can partly alleviate the regularity requirement on the solution.

Next we describe the deep Ritz method (DRM) \cite{yu2018deep}, which is based on the Ritz formulation. For problem \eqref{diri}, it is given by
\begin{equation}\label{deepritz}
\min_{u\in H^1(\Omega)} \tfrac{1}{2}(\kappa\nabla u,\nabla u)-f(\mathbf{x}_0)u(\mathbf{x}_0)+ \tfrac{\lambda}{2}\|u-h\|_{L^2(\partial\Omega)}^2.
\end{equation}
Here, using the definition of $\delta(\mathbf{x-x}_0)$, we eliminate the singularity in the source formally. By Sobolev embedding theorem \cite{Adams:2003}, the loss is well-posed in the one-dimensional case, but it is ill-defined when $d\geq2$, since the point evaluation $u(\mathbf{x}_0)$ is not defined for a function $u\in H^1(\Omega)$ when $d\geq2$.
Numerically, it does not perform well for problem \eqref{problem}, which is also confirmed by the numerical experiments below.

The weak adversarial network (WAN) \cite{Zang:2020} is based on the Galerkin weak formulation. For problem \eqref{diri}, it is derived as follows. Let $v\in H_0^1(\Omega)$ be any test function. Then by integration by parts,
\begin{equation}
	\int_{\Omega}(\kappa\nabla  u\cdot\nabla v-\delta(\mathbf{x}-\mathbf{x}_0) v-gv)\,{\rm d}\mathbf{x}=0.
\end{equation}
Thus, the weak solution $u$ solves the following min-max problem
\begin{equation}
	\min_{u\in H^1(\Omega): Tu = h}\max_{v\in H^1_0(\Omega)}\dfrac{\int_{\Omega}(\kappa\nabla u \cdot\nabla v-\delta (\mathbf{x}-\mathbf{x}_0)v-gv)\,{\rm d}\mathbf{x}}{\Big(\int_\Omega |\nabla v|^2{\rm d}\mathbf{x}\Big)^{1/2}}.
\end{equation}
Zang et al \cite{Zang:2020} proposed to employ two NNs for the trial function $u$ and test function $v$.
Due to its reliance on the weak formulation on $H^1(\Omega)\times H_0^1(\Omega)$, WAN suffers from the same
issue as DRM, since it is not well defined on the space $H^1(\Omega)\times H_0^1(\Omega)$.
Additionally, since the resulting min-max problem involves two optimization steps, the
training is often more delicate.
	
To address these challenges, several strategies have been proposed. One recent proposal is self-adaptive NNs \cite{mcclenny2020self}. In this approach, each term in the loss function is given a weight, and these weights are then adjusted to balance the effect of each term on the objective. This idea was recently employed to treat the influence brought by the singularity in PINN \cite{huang2021solving}, and the resulting method is called SAPINN. Specifically, one modifies the loss function into
\begin{equation*}
		L(u)=\lambda_1L_{r_1}(u)+\lambda_2L_{r_2}(u)+\lambda_3L_{b}(u),
\end{equation*}
where the three terms are respectively given by
\begin{align*}
L_{r_1}(u)&=\|-\nabla\cdot (\kappa\nabla u)-\delta_H(\mathbf{x}-\mathbf{x}_0)-g\|_{L^2(B(\mathbf{x}_0,H))}^2,\\
L_{r_2}(u)&=\|-\nabla\cdot (\kappa\nabla u)-g\|_{L^2(\Omega\backslash B(\mathbf{x}_0,H))}^2\quad \mbox{and}\quad
L_{b}(u)=\|u-h\|_{L^2(\partial\Omega)}^2,
\end{align*}
and $\lambda_i$, $i=1,2,3$, are positive weights to be adjusted, and $\delta_H(\mathbf{x}-\mathbf{x}_0)$ is a smooth approximation to the Dirac delta function $\delta(\mathbf{x}-\mathbf{x}_0)$ in a small ball $B(\mathbf{x}_0,H)$ (a ball of radius $H$, centered at $\mathbf{x}_0$). In this work we follow \cite{HOSSEINI2016423} and choose the approximation $\delta_H$ to be
$$ \delta_H(\mathbf{x})=
	\left\{
	\begin{aligned}
		\frac{12}{\pi H^d}\left(5\left(\frac{|\mathbf{x}|}{H}\right)^2-8\frac{|\mathbf{x}|}{H}+3\right),\quad|\mathbf{x}|\leq H,\\
		0,\quad|\mathbf{x}|>H.\\
	\end{aligned}\right. $$
The partition of the domain $\Omega$ into two disjoint parts with different penalty parameters allows better compensation for the singular behavior, in a manner similar to weighted Sobolev spaces. The experimental results in \cite{huang2021solving} indicate that this method can improve the accuracy of the solution by PINN to some extent, but the increase in the number of hyperparameters (i.e., $\lambda_i$s) also complicates the optimization process, due to the need for self adaptation during training, which in turn makes the overall training lengthier; see Section \ref{sec:numer} for a detailed comparison.

\section{Singularity Splitting Deep Ritz method}\label{sec:method}
In this section we propose a novel numerical method for solving the elliptic problem \eqref{problem}.
In practice, in addition to the Dirichlet boundary condition \eqref{diri}, equation \eqref{problem} may also be equipped with a Neumann or Robin one (possibly different conditions on different parts of the boundary $\partial\Omega$):
\begin{subequations}
\begin{align}
     \kappa\frac{\partial u}{\partial n}&=h,\quad \mbox{on }\partial\Omega,\label{neum}\\
        \alpha u+\kappa\frac{\partial u}{\partial n}&=h,\quad \mbox{on }\partial\Omega,
        \end{align}
    \end{subequations}
with the boundary data $h\in L^2(\partial\Omega)$ or its subspace, and $\alpha>0$. Below we
focus the discussion on the pure Dirichlet and Neumann boundary conditions, since the Robin
case can be treated similarly. In the Neumann case, one additionally needs to impose the
compatibility condition $\int_\Omega (f\delta + g){\rm d}\mathbf{x}+ \int_{\partial\Omega} h {\rm d}S = 0,$
so that a solution does exist. Then the solution is only unique up to an arbitrary constant,
and to have uniqueness, one can normalize the solution, e.g., $\int_\Omega u {\rm d}\mathbf{x}=0$
or $\int_{\partial\Omega} u{\rm d}S=0$, or just specify the value of the solution at one point
in the domain / on the boundary. All these considerations have to be incorporated into the loss function.

\subsection{Solution splitting}\label{ssec:split}
First we develop the singularity splitting technique for point and line sources.
One key challenge in the numerical solution of problem \eqref{problem} lies in
the low solution regularity, and thus a direct approximation of the solution $u$
via NNs can be inefficient. In order to overcome the challenge, we exploit
analytic properties of the solution $u$, or more specifically extract the solution
singularity directly and employ it in the solution procedure.  Note
that incorporating analytic insights into numerical procedures is a well
established computational paradigm, especially prominent in many analytic / semi-analytic
methods. Notable examples include locally refined meshes \cite{Apel:2011},
singular function enrichment \cite{Fix:1973}, singularity reconstruction
\cite{cai2001finite} and generalized FEM \cite{Fries:2010}. This work continues
along this long established tradition, and applies the idea in the context of neural PDE solvers for solving elliptic PDEs with singular sources.

First we describe the key analytic insights into the problem. Recall that the fundamental solution $ \Phi(\mathbf{x}) $ of the Laplace equation is given by \cite[p. 22]{evans10}
$$\Phi(\mathbf{x})=\left\{
	\begin{aligned}
		-\frac{1}{2\pi}\ln|\mathbf{x}|,&\quad d=2,\\
		\frac{1}{d(d-2)\alpha(d)|\mathbf{x}|^{d-2}},&\quad d\geq3,\\
	\end{aligned}\right.$$
where $ \alpha(d) =\frac{\pi^\frac{d}{2}}{\Gamma(\frac{d}{2}+1)}$ denotes the volume of the unit ball in $\mathbb{R}^d$. Note that $\Phi$ is smooth away from the singularity at the origin. By definition, $ \Phi(\mathbf{x}) $ satisfies
\begin{equation*}
     -\Delta \Phi(\mathbf{x-x}_0)=\delta(\mathbf{x-x}_0),\quad \mbox{in } \mathbb{R}^d\setminus\{\mathbf{x}_0\},
\end{equation*}
where $\delta(\mathbf{x}-\mathbf{x}_0)$ denotes the Dirac delta function concentrated at the point $\mathbf{x}_0\in\mathbb{R}^d$.

\subsubsection{Singularity reconstruction for point sources}

Using the fundamental solution $\Phi(\mathbf{x})$, we can split the solution $u $ to problem
\eqref{problem} with one point source $S(\mathbf{x})=c_0\delta(\mathbf{x-x}_0)$ into
\begin{equation}\label{split}
    	u=c_0\kappa^{-1}\Phi+v,
\end{equation}
where $\Phi$ is an abbreviation of $\Phi(\mathbf{x-x}_0)$. The extension of the approach and
the discussion below to the case of multiple Dirac functions is straightforward. In the splitting
\eqref{split}, the term $c_0\kappa^{-1}\Phi$ captures the leading singularity of the solution
$u$ near the point source $\mathbf{x}_0$ and $ v $ is the regular part. Since the fundamental solution
$\Phi(\mathbf{x-x}_0)\neq0$ when $\mathbf{x}\in\partial\Omega$, the boundary condition and the
source term have to be modified accordingly. For example, for the Dirichlet problem, by
substituting the splitting \eqref{split} into \eqref{problem}, we get the
following governing equation for the regular part $v$:
    \begin{equation}\label{prob}
    	\left\{
    	\begin{aligned}
    		-\nabla\cdot (\kappa\nabla v)&=F,&& \mbox{in }\Omega,\\
    		v&=h-c_0\kappa^{-1}\Phi,&& \mbox{on }\partial\Omega,\\
    	\end{aligned}
    	\right.
    \end{equation}
	with the modified source $F$ given by
	\begin{equation}\label{right}
		F=g-c_0\kappa^{-1}[\Delta\kappa \Phi+\nabla\kappa\cdot \nabla\Phi]+c_0\kappa^{-2}|\nabla\kappa|^2\Phi.
\end{equation}

Next, we state the Ritz variational formulation for the regular part $v$, on which DRM is based. Let $ \tilde h= h -c_0\kappa^{-1}\Phi$ in the Dirichlet case, and $\tilde h = h - c_0\kappa\partial_n(\kappa^{-1}\Phi)$ in the Neumann case, be the modified boundary data. Then the Ritz variational formulation is posed on $H^1(\Omega)$ and given by
\begin{equation}\label{eqn:Ritz-regular}
   L(v)=\left\{\begin{aligned}\tfrac{1}{2}(\kappa\nabla v,\nabla v)-(F,v), & \quad \mbox{Dirichlet},\\
   \tfrac{1}{2}(\kappa\nabla v,\nabla v)-(F,v)-(\tilde h,v)_{L^2(\partial\Omega)}, & \quad \mbox{Neumann}.
   \end{aligned}\right.
\end{equation}

Let $T$ be the trace operator. Then we have the following well known result \cite{evans10}.
\begin{proposition}
If $F\in (H^1(\Omega))'$ and $\tilde h\in H^\frac12(\partial\Omega)$ and $\tilde h\in (H^\frac12(\partial\Omega))'$ in the Dirichlet and Neumann
case, respectively. Then the regular part
$ v$ is a minimizer of the functional $L(v)$ defined in \eqref{eqn:Ritz-regular} over the
set $\{v\in H^1(\Omega): Tv=\tilde h \}$ and $H^1(\Omega)$, in the Dirichlet and Neumann
case, respectively.
\end{proposition}

The regularity assumptions on $F$ and $\tilde h$ will be analyzed in Section \ref{ssec:reg} below.
Note that in the Neumann case, a solution $u$ exists only under suitable compatibility conditions
and to achieve the uniqueness, one may specify the value of the solution $u$ at a selected point
$x^*\in\overline{\Omega}$, e.g., $u({x^*})=a$. Then it has to be supplied with the condition
$v({x^*})=\tilde a \equiv a-\kappa({x^*})^{-1}c_0\Phi({x^*})$, which is then used as a
constraint in minimizing the functional $L(v)$.

\subsubsection{Well-posedness of the modified problem \eqref{prob}}\label{ssec:reg}
Note that in the presence of a point source, the solution $u$ contains the singularity of the fundamental
solution $\Phi$, which belongs only to $W^{1,p}(\Omega)$ for $p<\frac{d}{d-1}$ \cite[Theorem 1.1, p. 3]{GruterWidman:1982},
and as the dimension $d$ grows, the Sobolev regularity of $\Phi$ deteriorates.
In particular, it does not have the $H^1(\Omega)$ regularity when $d\geq2$, and thus the standard Ritz
type formulation fails to make sense. Next, we discuss the smoothness of the regular part $v$.
This depends on the regularity of $\kappa$ and $g$ and that of the coefficient $\kappa$. Below we discuss
the regularity of the modified source $F$ and the regular part $v$ for one single point source. The next
result summarizes the regularity of $F$.  The special case of $\kappa$
being locally constant in the neighborhood of the singularities is common in practice,
e.g., the standard Poisson equation, which is arguably the most widely studied case for
elliptic PDEs with point sources \cite{Babuska:1971,Scott:1974,KoppWohlmuth:2014}. The analysis
also sheds interesting insides into regularity restriction and the impact of the dimension $d$.
\begin{proposition}\label{prop:reg-F}
Fix $d\geq2$. Let $\kappa\in W^{2,r}(\Omega)$ with $r\geq 2$, and $g\in L^2(\Omega)$. Then there holds
\begin{equation*}
    F\in L^s(\Omega), \quad \mbox{with } \left\{\begin{aligned}
     s<\frac{d}{d-1},  &\quad r\geq d, ~d \geq 2,\\
    s<\min\Big(\frac{dr}{(d -2)r +d},\frac{dr}{(d-2)r+2(d-r)}\Big) &\quad r<d,~ d \geq 3.
    \end{aligned}\right.
\end{equation*}
Further, if $\kappa$ is locally constant in the neighbourhood of the singularity at $\mathbf{x}_0$, then $F\in L^2(\Omega)$.
\end{proposition}
\begin{proof}
The assertion follows from direct computation with the standard Sobolev embedding theorem \cite[Theorem 4.12, p. 85]{Adams:2003}. It is known that $\Phi\in W^{1,p}(\Omega)$ for any $p<\frac{d}{d-1}$ \cite[Theorem 1.1, p. 3]{GruterWidman:1982}. Meanwhile, by Sobolev embedding theorem,
\begin{equation}\label{eqn:embedding}
    W^{1,r}(\Omega)\hookrightarrow L^{r^*}(\Omega),\quad \mbox{with }
    \left\{\begin{aligned}
     r^*=\infty, &\quad  r>d,\\
     r^*<\infty, &\quad r=d,\\
     r^*=\frac{dr}{d-r}, & \quad r<d.
\end{aligned} \right.
\end{equation}
Using the embedding \eqref{eqn:embedding} and the regularity $ \Phi \in W^{1,p}(\Omega)$
with $p<\frac{d}{d-1}$ \cite[Theorem 1.1, p.3]{GruterWidman:1982}, we have
\begin{equation}\label{eqn:embedd-Phi}
   \Phi \in L^{p^*}(\Omega), \quad \mbox{with } \left\{\begin{aligned}
    p^*<\infty, &\quad   d=2, \\
    p^*<\frac{d}{d-2}, &\quad d\geq 3.
   \end{aligned}\right.
\end{equation}
In view of the expression \eqref{right} of $F$, it suffices
to analyze the three terms $\Phi\Delta \kappa$, $\nabla \kappa \cdot \nabla \Phi$ and
$|\nabla \kappa|^2\Phi$. By H\"{o}lder's inequality, we deduce
\begin{equation*}
    \Phi\Delta \kappa \in L^s(\Omega), \quad \mbox{with } \left\{\begin{aligned}
        s <r, & \quad d=2,\\
        s<\frac{dr}{(d-2)r+d}, &\quad d\geq 3.
    \end{aligned}\right.
\end{equation*}
Similarly, for the term $\nabla \kappa\cdot \nabla \Phi$ (noting the condition $r\geq2$), we have
\begin{equation*}
    \nabla \kappa\cdot \nabla \Phi \in L^s(\Omega), \quad \mbox{with } \left\{\begin{aligned}
     s<\frac{d}{d-1},  &\quad r\geq d, ~d \geq 2,\\
    s<\frac{dr}{(d -2)r +d}, &\quad r<d,~ d \geq 3.
    \end{aligned}\right.
\end{equation*}
Last, it follows from the embeddings \eqref{eqn:embedding} and \eqref{eqn:embedd-Phi},
and the generalized H\"{o}lder inequality that
\begin{equation*}
  |\nabla\kappa|^2\Phi \in L^q(\Omega),\quad \mbox{with }\left\{\begin{aligned}
    q< \infty, &\quad d = 2,\\
    q<\frac{d}{d-2}, &\quad r\geq d,~ d\geq 3,\\
    q<\frac{dr}{(d-2)r+2(d-r)}, &\quad r<d,~ d\geq 3.\\
  \end{aligned}\right.
\end{equation*}
Now combining the preceding three results yields the first assertion. The remaining assertion follows from direct computation.
\end{proof}

Based on the regularity result of the modified source $F$ in Proposition \ref{prop:reg-F}, we can
discuss the well-posedness of the Ritz variational formulation of problem \eqref{prob}, which lays
the foundation for applying DRM. For a general variable coefficient $\kappa$, when $d=2,3$, the formulation is indeed well posed. Then by the
standard elliptic regularity theory, we have $v\in H^1(\Omega)$, whereas the regularity of the
solution $u$ to problem \eqref{problem} with a point source is inherently limited to $u\in W^{1,p}(\Omega)$ for any
$p<\frac{d}{d-1}$. One can also deduce the precise Sobolev regularity of the regular part $v$:
in the 2D case $v$ belongs to $W^{2,s}(\Omega)$ for any $s<2$, just falling short of
$H^2(\Omega)$, whereas in the 3D case, $v $ can at most have a regularity $W^{2,\frac{3}{2}}(\Omega)$
(when $r=\infty$). This shows the benefit of singularity splitting for problem
\eqref{problem}. When $d\geq4$, it may be ill-defined, except the special case
$\kappa$ being locally constant. In the latter case, $v$ does belong to $H^2(\Omega)$ for any $d\geq2$, and thus
the regular part $v$ indeed enjoys much better regularity than the solution $u$.
\begin{corollary}\label{cor:well-posed}
Let $\kappa\in W^{2,r}({\Omega})$ with $r\geq d$, $h\in H^\frac{1}{2}(\partial\Omega)$ and $g\in L^2(\Omega)$.  Then
problem \eqref{prob} has a weak solution $v\in H^1(\Omega)$ for $d=2,3$. Further, if $\kappa$ is
locally constant in a neighborhood of $\mathbf{x}_0$, then
it has a weak solution $v\in H^1(\Omega)$ for any $d\geq2$.
\end{corollary}
\begin{proof}
First, since the singularity point $ \mathbf{x}_0\in \Omega$, we have $\kappa^{-1}\phi\in H^\frac{1}{2}(\partial\Omega)$, and hence $\tilde h \in H^\frac12(\partial\Omega)$. By Sobolev embedding \eqref{eqn:embedding}, we have
\begin{equation*}
   H^1(\Omega)\hookrightarrow L^s(\Omega), \quad \mbox{with } \left\{\begin{aligned}
       s <\infty, & \quad d=2,\\
       s = \frac{2d}{d-2}, &\quad d\geq 3.
   \end{aligned}\right.
\end{equation*}
Therefore, by duality, we have
\begin{equation*}
    L^{s^*}(\Omega)\hookrightarrow (H^1(\Omega))', \quad \mbox{with }\left\{\begin{aligned}
       s^* >1, & \quad d=2,\\
       s^* = \frac{2d}{d+2}, &\quad d\geq 3.
   \end{aligned}\right.
\end{equation*}
Then for $d=2, r\geq2$, clearly $F\in L^s(\Omega)$, for any $s<2$. Similarly, for
$d=3, r>2$, we have $ F\in L^s(\Omega)$, with $s<\frac{3r}{r+3}$, and $\frac{3r}{r+3}>\frac{2d}{d+2}$.
Meanwhile, for $d=4, r>4$, we have $F\in L^s(\Omega)$ for any $s<\frac{4}{3}=
\frac{2d}{d+2}$. Thus, for sufficiently smooth $g$, we have $F\in (H^1(\Omega))'$
for $d=2,3$ only. Meanwhile, if $\kappa$ is locally constant in a neighborhood of
the singularity at $\mathbf{x}_0$, then $F\in L^2(\Omega)\hookrightarrow (H^1(\Omega))'$,
for any $d\geq2$.
\end{proof}

\subsubsection{Singularity splitting for line sources}

Next we discuss the case of line sources. We describe the procedure for a line source supported
on a line segment in $\mathbb{R}^d$, $d\ge3$.
The construction below follows the procedure in \cite{gjerde2019splitting}, which also contains a detailed construction of the singularity $\Phi_L$ in $\mathbb{R}^3$.
Nonetheless, to the best of our knowledge, the resulting expressions appear to be new, since the
numerical approximation of high-dimensional PDEs has not received the due attention.
Specifically, we parameterize the line segment $\Lambda \subset\Omega$ connecting two points
$\mathbf{a},\mathbf{b}\in \Omega$ by $\mathbf{y}=\mathbf{a}+ t\boldsymbol{\tau}$ for $t\in(0,L)$,
with $L=|\mathbf{b-a}|$ and $\boldsymbol{\tau}=(\mathbf{b-a})/L$ the unit tangent vector of
$\Lambda$. Consider the elliptic  problem \eqref{problem}
with $\delta$ being replaced with $\delta_\Lambda$. Then we seek a function $\Phi_L$ so that
\begin{equation*}
    -\int_\Omega \Delta \Phi_L v\,{\rm d} \mathbf{x} = \int_\Lambda v\,\mathrm{d}s,\quad \forall v\in C(\overline{\Omega}).
\end{equation*}
A natural candidate for $\Phi_L$ is the convolution of $\delta_\Lambda$ and $\Phi$, i.e.,
\begin{equation*}
\Phi_L(\mathbf{x}) = \int_\Omega \delta_\Lambda(\mathbf{y})\Phi(\mathbf{x-y})\,{\rm d}\mathbf{y}=\int_0^L\frac{1}{d(d-2)\alpha(d)}|\mathbf{x}-(\mathbf{a}+\boldsymbol{\tau}t)|^{-(d-2)}\,{\rm d}t.
\end{equation*}
So it suffices to evaluate the integral ${\rm I}:=\int_{0}^{L}|\mathbf{x}-(\mathbf{a}+\boldsymbol{\tau}t)|^{2-d}{\rm d}t $.
Note that in $\mathbb{R}^3$, with $r_a=|\mathbf{x-a}|$ and $r_b=|\mathbf{x-b}|$, we can compute directly \cite[p. 1722]{gjerde2019splitting}
\begin{equation*}
   {\rm I} =\ln\frac{r_b+L+\boldsymbol{\tau}\cdot(\mathbf{a-x})}{r_a+\boldsymbol{\tau}\cdot(\mathbf{a-x})}.
\end{equation*}
For $d\geq4$, by first setting $\alpha^2=r_a^2-((\mathbf{a-x})\cdot\boldsymbol{\tau})^2$, and then substituting $t=\alpha \tan s$, we have
\begin{align*}
	{\rm I}&=\int_{0}^{L}[t^2+2(\mathbf{a-x})\cdot\boldsymbol{\tau}t+r_a^2]^{1-\frac{d}{2}}{\rm d}t=\int_0^L(t+(\mathbf{a-x})\cdot\tau)^2-((\mathbf{a-x})\cdot\tau)^2+r_a^2{\rm d}t\\
  &=\int_{(\mathbf{a-x})\cdot\boldsymbol{\tau}}^{(\mathbf{a-x})\cdot\boldsymbol{\tau}+L}\!\!\!{(t^2+\alpha^2)^{1-\frac{d}{2}}}{\rm d}t=
 \dfrac{1}{\alpha^{d-3}}\int_{\arctan\frac{(\mathbf{a-x})\cdot\boldsymbol{\tau}}{\alpha}}^{\arctan\frac{(\mathbf{a-x})\cdot\boldsymbol{\tau}+L}{\alpha}}\!\!\!\cos^{d-4}s{\rm d}s.
\end{align*}
Then using the following integral identities \cite[2.513, p. 153]{GradshteynRyzhik:2015}
\begin{align*}
  \int \cos ^{2n}x{\rm d}x &= \frac{1}{2^{2n}}{2n \choose n}x + \frac{1}{2^{2n-1}}\sum_{k=0}^{n-1}{2n \choose k}\frac{\sin(2n-2k)x}{2n-2k},\\
  \int \cos ^{2n+1}x{\rm d}x &= \frac{1}{2^{2n}}\sum_{k=0}^{n}{2n+1 \choose k}\frac{\sin(2n-2k+1)x}{2n-2k+1},
\end{align*}
we obtain
 \begin{align*}
 {\rm I}=\frac{1}{\alpha^{d-3}}
	    \left\{
		\begin{aligned}
		&\arctan\frac{(\mathbf{a-x})\cdot\boldsymbol{\tau}+L}{\alpha}-\arctan\frac{(\mathbf{a-x})\cdot\boldsymbol{\tau}}{\alpha},\quad d=4,\\
		&\dfrac{1}{2^{d-5}}\sum_{i=0}^{m-3}{d-4\choose i}\frac{\sin[(d-4-2i)\arctan\frac{(\mathbf{a-x})\cdot\boldsymbol{\tau}+L}{\alpha}]-\sin[(d-4-2i)\arctan\frac{(\mathbf{a-x})\cdot\boldsymbol{\tau}}{\alpha}]}{d-4-2i}\\
		&\quad+\dfrac{1}{2^{d-4}}{d-4\choose m-2}\left[\arctan\frac{(\mathbf{a-x})\cdot\boldsymbol{\tau}+L}{\alpha}-\arctan\frac{(\mathbf{a-x})\cdot\boldsymbol{\tau}}{\alpha}\right],\quad d=2m\geq6,\\
		&\frac{1}{2^{d-5}}\sum_{i=1}^{m-2}{d-4\choose i}\frac{\sin[(d-4-2i)\arctan\frac{(\mathbf{a-x})\cdot\boldsymbol{\tau}+L}{\alpha}]-\sin[(d-4-2i)\arctan\frac{(\mathbf{a-x})\cdot\boldsymbol{\tau}}{\alpha}]}{d-4-2i},\\&\quad d=2m+1\geq5.
	\end{aligned}
     \right.
\end{align*}
Combining these identities yields an explicit expression for $\Phi_L(\mathbf{x})$.
Using the function $\Phi_L(\mathbf{x})$, we can split the solution $u $ to problem \eqref{problem} into
\begin{equation}\label{split-line}
    	u=\kappa^{-1}f\Phi_L+v.
\end{equation}
In the splitting \eqref{split-line}, the term $\kappa^{-1}f\Phi_L$ captures the leading singularity of the solution $u$ near the line sources and $ v $ is the regular part. Since the fundamental solution $\Phi_L(\mathbf{x})\neq0$ when $\mathbf{x}\in\partial\Omega$, the boundary condition and the source term have to be modified accordingly. For example, for the Dirichlet problem, by substituting the splitting \eqref{split-line} into \eqref{problem} we can get the following governing equation for the regular part $v$:
    \begin{equation}\label{prob-line}
    	\left\{
    	\begin{aligned}
    		-\nabla\cdot (\kappa\nabla v)&=F,&&\quad \mbox{in }\Omega,\\
    		v&=h-\kappa^{-1}f\Phi_L,&&\quad \mbox{on }\partial\Omega,\\
    	\end{aligned}
    	\right.
    \end{equation}
	with the modified source $F$ given by
	\begin{equation}\label{right-line}
		F=g+\Delta f\Phi_L+2\nabla f\cdot \nabla\Phi_L-\kappa^{-1}[\Delta\kappa f\Phi_L+f\nabla\kappa\cdot \nabla\Phi_L +\Phi_L\nabla\kappa\cdot\nabla f]+\kappa^{-2}|\nabla\kappa|^2 f\Phi.
\end{equation}

Note that the regularity of the regular part $v$ to problem \eqref{prob-line} can also be
analyzed similarly as in Section \ref{ssec:reg}, under suitable regularity assumption
on the density $f$. Indeed, line sources can be obtained by taking Laplacian of a lower
dimensional fundamental solution in the high-dimensional ambient space. Similar to Corollary \ref{cor:well-posed}, generally there is still a restriction on the dimension $d$ (i.e., $d=2,3,4$ only) in the presence of line sources, in order for the variational problem to be well-posed, but again the restriction disappears when $f$ and $\kappa$ are locally constant in the neighborhood of singularity support, i.e., the line segment. See also the work
\cite{gjerde2019splitting} for some insightful discussions on weighted Sobolev regularity
in the 3D case. We leave a precise regularity analysis, including in weighted
Sobolev spaces, to a future work. Using the convolution strategy, one can also derive the
expression for more general singularities, e.g., singularities on squares / cubes in
high-dimensional space. Due to the complexity of the resulting expressions, we do not
further pursue the issue below.

\subsection{Deep Ritz method}
Now we describe the deep Ritz method (DRM) \cite{yu2018deep} for approximating the regular part $v$ and the details about the practical implementation. First we discuss the Dirichlet case, and will comment on the Neumann case at the end, which only requires minor changes.
The energy functional of problem \eqref{prob} requires minimizing the functional $L(v)$ in \eqref{eqn:Ritz-regular} over the set $\{v\in H^1(\Omega): Tv = \tilde h\}$. However, due to the global nature of NNs, it is highly nontrivial to construct an NN function that satisfies the Dirichlet boundary condition $Tv=\tilde h$ exactly, although theoretically this is highly desirable for improving the error estimates of the method \cite{MullerZeinhofer:2021exact}. In order to enforce the Dirichlet boundary condition $Tv=\tilde h$, we adopt the standard penalty formulation
\begin{equation}
    L_\sigma(v)=\tfrac{1}{2}(\kappa \nabla v,\nabla v) - (F, v) +\tfrac{\sigma}{2}\|v-\tilde h\|_{L^2(\partial\Omega)}^2,\quad v\in H^1(\Omega),
\end{equation}
where $\sigma>0$ is the penalty parameter. This choice is standard in
the current practice of NN based PDE solvers \cite{RAISSI2019686,yu2018deep}.
Intuitively, as the penalty parameter $\sigma$ increases to infinity, the equality
constraint representing the Dirichlet boundary condition is better satisfied; see
Theorem \ref{thm:penalization} and the remark thereafter for precise statements.

Now we describe NNs for approximating the regular part $v$, which is the central component of DRM. We employ standard fully connected feed forward NNs, i.e., functions $ f_\theta: \mathbb{R}^d\rightarrow \mathbb{R}$, with the network parameter $\theta$. Specifically, it is defined recursively by
\begin{align*}
   \mathbf{x}_0&=\mathbf{x},\\
 \mathbf{x}_\ell&=\rho(\mathbf{A}_\ell\mathbf{x}_{\ell-1}+\mathbf{b}_\ell),\quad \ell=1,2,\cdots,L-1,\\
	f_\theta(\mathbf{x})&=\mathbf{A}_L\mathbf{x}_{L-1}+\mathbf{b}_L,
\end{align*}
where $\mathbf{A}_\ell\in\mathbb{R}^{n_\ell\times n_{\ell-1}}, \mathbf{b}_\ell\in\mathbb{R}^{n_\ell},\, \ell=1,2,\cdots,L$. Clearly we fix $n_0=d$ and $n_L=1$. In the construction, the nonlinear function $\rho:\mathbb{R}\to \mathbb{R}$ is called the activation function, and it is applied componentwise to a vector. $L$ is called the depth, and $W:=\max\{n_\ell,\ell=0,1,\cdots,L\}$ is called the width of the network. The set of parameters of the neural network, i.e., $\mathbf{A}_\ell,\,\mathbf{b}_\ell\,\,(\ell=1,2,\cdots,L)$, are trainable and are often stacked into a large vector $\theta$. The training is carried out by minimizing a suitable loss function, which will be specified below. In the construction, $\mathbf{x}_0$ is called the input layer, $\mathbf{x}_\ell$, $\ell=1,2,\cdots,L-1$, are called the hidden layer and $f_\theta(\mathbf{x})$ is the output layer, where the subscript $\theta$ explicitly indicates the dependence of the NN on the parameter vector $\theta$. The total number of parameters in the neural network $f_\theta$ is
$\sum_{\ell=1}^Ln_\ell(n_{\ell-1}+1).$

For the activation function $\rho$, there are many possible choices. The most frequently used one in computer vision is the rectified linear unit (ReLU), defined by $\rho(x)=\max(x,0)$  \cite{LeCunBottou:1998}. However, it is not smooth enough for the use in (SS)DRM, since the formulation requires twice differentiability of the activation function $\rho$: one spatial derivative, and the optimizer requires another derivative in the NN parameter $\theta$. In the context of neural PDE solvers, hyperbolic tangent $\rho (x)=\frac{{\rm e}^x-{\rm e}^{-x}}{{\rm e}^x+{\rm e}^{-x}}$ and logistic $ \rho(x)=\frac{1}{1+{\rm e}^{-x}} $ are frequently used \cite{RAISSI2019686,CuomoSchiano:2022}. In the numerical experiments, we employ the hyperbolic tangent as the activation function. Note that under fairly generous conditions on $\rho$, the set of NN functions can approximate any $L^p(\Omega)$ function arbitrarily well, as the width $W$ and / or depth $L$ tend to infinity, a property commonly known as universal approximation property. Furthermore, quantitative estimates on the approximation error in various norms have also been derived in recent years \cite{yarotsky2017error,I2020Approximation}.  Below we denote the collection of NN functions of depth $L$, with the total number of nonzero parameters $N_\theta$, and each the parameter being bound by $R$, with the activation function $\rho$ by $\mathcal{N}_\rho(L,N_\theta,R)$, i.e.,
\begin{equation}\label{eqn:nn-set}
  \mathcal{N}_\rho(L,N_\theta,R)   = \{v_\theta: v_\theta ~ \mbox{has a depth } L,\, |\theta|_0\leq N_\theta, \, |\theta|_{\ell^\infty}\leq R\},
\end{equation}
where $|\cdot|_{\ell^0}$ and $|\cdot|_{\ell^\infty}$ denote the number of nonzero entries in and the maximum norm of a vector, respectively.
Below we also use the shorthand notation $\mathcal{A}$ to denote this collection of functions (with a fixed architecture).

In a neural PDE solver, we use an element $v_\theta$ from the set $\mathcal{A}$ to approximate the minimizer $v_\sigma\in H^1(\Omega)$ of the functional $L_\sigma(v)$, i.e.,
\begin{equation*}
  \theta^* = \arg \min_{\theta} L_\sigma(v_\theta),
\end{equation*}
and then set $v_{\theta^*}\in\mathcal{A}$ as an approximation to $v_\sigma$. Since the parameter $\theta$
is finite-dimensional and lives on compact set (due to the $\ell^\infty$ bound $|\theta|_{\ell^\infty}
\leq R$), for a smooth activation function $\rho$, the loss $L_\sigma(v_\theta)$ is continuous in $\theta$,
which directly implies the existence of a global minimizer $\theta^*$. Note that the loss $L_\sigma(v_\theta) $
involves potentially high-dimensional integrals (i.e., in $\Omega$ and on $\partial\Omega$), and thus
direct computation is likely still intractable. Instead, in practical computation, we approximate the
integrals using quadrature, often by means of Monte Carlo, especially in the high-dimensional case.
Specifically, we rewrite the loss $L_\sigma$ as
\begin{equation}\label{eqn:loss-cont}L_\sigma(v_\theta)=|\Omega|\mathbb{E}_{X\sim U(\Omega)}\Big[\frac{\kappa(X)|\nabla v_\theta(X)|^2}{2}-F(X)v_\theta(X)\Big] + \frac{\sigma |\partial\Omega|}{2}\mathbb{E}_{Y\sim U(\partial\Omega)}\big[|v_\theta(Y)-\tilde h(Y)|^2\big],
\end{equation}
where $U(\Omega)$ and $U(\partial\Omega)$ denote the uniform distributions on the domain $\Omega$ and the boundary $\partial\Omega$, respectively, and $\mathbb{E}_{U(\Omega)}$ denotes taking expectation with respect to the uniform distribution $U(\Omega)$. This consideration leads to the following empirical loss
\begin{equation}\label{eqn:loss-emp}\widehat{L}_\sigma(v_\theta)=\dfrac{|\Omega|}{N_r}\sum_{i=1}^{N_r}\left[\dfrac{\kappa(X_i)|\nabla v_\theta(X_i)|^2}{2}-F(X_i)v_\theta(X_i)\right]+\dfrac{|\partial\Omega|}{N_{b}}\frac{\sigma}{2}\sum_{j=1}^{N_{b}}(v_\theta(Y_j)-\tilde h(Y_j))^2,
\end{equation}
where the sampling points $\{X_i\}_{i=1}^{N_r}$ and $\{Y_j\}_{j=1}^{N_{b}}$ are identically and independently distributed (i.i.d.) random variables, distributed uniformly on the domain $\Omega$ and the boundary $\partial\Omega$, respectively, i.e., $\{X_i\}_{i=1}^{N_r}\sim U(\Omega)$ and $\{Y_j\}_{j=1}^{N_{b}}\sim U(\partial\Omega)$. Note that the distributions are taken to be uniform according to the continuous loss \eqref{eqn:loss-cont}. In practice, it is also possible to use nonuniform distributions, which corresponds to solving the PDEs in suitable weighted spaces and can be beneficial for certain problems; see \cite{WuZhuLu:2022,TangWanYang:2022} for numerical studies along this direction.
Note that the resulting optimization problem in $\theta$ reads
\begin{equation}\label{eqn:min-emp}
    \widehat{\theta^*} = \arg\min_\theta \widehat{L}_\sigma(v_\theta)
\end{equation}
with $v_{\widehat{\theta^*}}\in\mathcal{A}$ being the NN approximation. The continuous loss
$ L_\sigma(v_\theta) $ and empirical loss $ \widehat{L}_\sigma(v_\theta)$ have different minimizers,
due to the introduction of additional quadrature errors. We shall discuss this issue briefly in
Section \ref{sec:theory}.

The final optimization problem \eqref{eqn:min-emp} is highly nonconvex, due to the nonlinearity
of the NN $v_\theta$ in the NN parameter $\theta$. Therefore, it is very challenging to find a global minimizer of the loss $\widehat{L}_\sigma(v_\theta)$. Nonetheless, in practice, simple optimization
algorithms, e.g., stochastic gradient descent and its variants (e.g. ADAM and
AdaGrad) \cite{BottouCurtis:2018}, seem to work fairly well, as testified by the great empirical
successes \cite{CuomoSchiano:2022}. In the numerical experiments, we will employ the limited memory
BFGS algorithm \cite{ByrdLu:1995} to minimize the empirical loss $\widehat{L}_\sigma(v_\theta)$,
which has been a very popular choice for NN based PDE solvers \cite{RAISSI2019686,CuomoSchiano:2022}.

In the Neumann case, the empirical loss $\widehat{L}_\sigma(v_\theta)$ takes a similar form:
\begin{equation*}
  \widehat{L}_\sigma(v_\theta) = \frac{|\Omega|}{N_r}\sum_{i=1}^{N_r}\Big[\frac12\kappa(X_i)|\nabla v_\theta(X_i)|^2 - F(X_i)v_\theta(X_i)\Big] - \frac{|\partial\Omega|}{N_b}\sum_{j=1}^{N_b}\tilde h(Y_j) v_\theta(Y_j)+ \frac{\sigma}{2}|v_\theta(\mathbf{x}^*)-\tilde a|^2,
\end{equation*}
where the last term is to approximately enforce the point value evaluation $v(\mathbf{x}^*)=
\tilde a$ (with $\mathbf{x}^*\in\overline{\Omega}$).

\begin{algorithm}[H]
\caption{The penalty method for SSDRM}
\begin{algorithmic}
  \State Set $ \sigma_1>0$, increasing factor $\eta>1$
\While{Stopping condition not met}
		\State Minimize the loss function $\widehat{L}_{\sigma_k}(v_\theta)$, initialized to $\widehat{\theta}_{k-1}^*$ and find the optimal $\widehat{\theta}_k$
   \State $ \sigma_{k+1} \leftarrow\eta\sigma_k$, and $ k \leftarrow k+1$
			\EndWhile
		\end{algorithmic}
	\end{algorithm}

In practice, the choice of the penalty parameter $\sigma$ in the objective is important since it
balances several competing effects: the larger is the value of $\sigma$, the smaller is the
error due to penalization, but the larger is the approximation error on the boundary, as indicated
by Theorem \ref{thm:error}. Thus, a good tradeoff between these different sources of errors is needed in order to
deliver high-quality NN approximations to the solution $u$ of the PDE \eqref{problem}. We have
followed a simple path following procedure, a well established procedure in optimization
\cite{AllgowerGeorg:2003,JiaoJinLu:2015}: we start with a small value $\sigma_1$, and then
after each loop update $\sigma$ geometrically: $\sigma_{k+1}=\eta\sigma_k$, for some fixed
$\eta>1$. By updating $\sigma_k$, the minimizer $\widehat{\theta}^*_k$ of the loss $\widehat{L}_{\sigma_k}(v_\theta)$
also approaches the solution of problem \eqref{prob}, and the path-following strategy allows
enforcing this progressively.
The overall procedure is shown in Algorithm 1. The use of the path following strategy requires the training for
multiple $\sigma$ values, which can potentially be time consuming, if done naively. Fortunately, the parameter $\theta$ of the $\sigma_{k+1}$-problem (i.e., minimizing $\widehat{L}_{\sigma_{k+1}}(v_\theta)$)  can be initialized to the converged parameter $\widehat{\theta}_k^*$ of $\sigma_k$-problem in order to warm start the optimization process.
This ensures that for each fixed $\sigma_k$ (except $\sigma_1$), the initial parameter configuration is close to the optimal one, so the
training loop only requires relatively few iterations to reach convergence.  This is also confirmed by the numerical experiments below.
	
\section{Numerical experiments and discussions}\label{sec:numer}
In this section, we present several numerical experiments to illustrate the performance of the
proposed SSDRM. Although it works for all three types of boundary conditions, we
present results mostly for the Dirichlet problem, and we use $f$ to denote the singularity strength for point sources / density for line sources (and thus omit $c_i$ for singularity strength) and $\mathbf{x}=(x_1,\ldots,x_d)\in\mathbb{R}^d$.
In the training, $N_r=10,000$ points
in the domain $\Omega$ and $N_b=400$ points on the boundary $\partial\Omega$ are randomly selected to form the empirical loss $\widehat{L}_\sigma(v_\theta)$, unless otherwise specified. In the path-following strategy for determining the penalty parameter $\sigma$, we take the initial penalty parameter $\sigma_1=20$,
and an increasing factor $\eta=1.5$. All the numerical experiments were carried out on a personal laptop (operating-system: Windows 10, with RAM 8.0GB, Intel(R) Core(TM) i7-10510U CPU, 2.3 GHz), with Python 3.9.7, with the popular software framework PyTorch. The gradient of the NN output with respect to the input $\mathbf{x}$ (i.e., spatial derivatives) and that of the empirical loss with respect to the NN parameter vector $\theta$ are computed via automatic differentiation, as implemented by  \texttt{torch.autograd}. The empirical loss $\widehat{L}_\sigma(v_\theta)$ is minimized using the off-the-shelf optimizer limited memory BFGS \cite{ByrdLu:1995}, as implemented in the SciPy library, with the default setting (tolerance =1.0e-9, no box constraint) and a maximum of 2500 iterations. To measure the accuracy of an approximation $\hat v$ of the regular part $v^*$, we use the standard relative $L^2(\Omega)$-error $e$ defined by
$e=\|v^*-\hat v\|_{L^2(\Omega)}/\|v^*\|_{L^2(\Omega)}$. The stopping condition is set to $e<\text{1.00e-3}$ and $\sigma_k\leq 1200$. The first
condition ensures that the approximation can achieve the desired accuracy and the second one is to terminate
the iteration after a finite number of loops.  The PyTorch source code for reproducing all the numerical experiments will be made available at the github link \url{https://github.com/hhjc-web/SSDRM.git}.

First, we showcase the approach on the case with one point source.
\begin{example}\label{exam:2d-1ps-Dirichlet}
The domain $\Omega$ is $\Omega=(-1,1)^2$, diffusion coefficient $\kappa=|\mathbf{x}|^2+1$, a point source at $\mathbf{x}_0=(0,0)$ with strength $f=2(|\mathbf{x}|^2+1)$,
source $g=(|\mathbf{x}|^2+1)\left(2x_2\sin x_1{\rm e}^{x_1x_2}+(1-x_1^2+x_2^2){\rm e}^{x_1x_2}\cos x_1\right)+\frac{2}{\pi}-4x_1x_2\cos x_1{\rm e}^{x_1x_2}+2x_1\sin x_1{\rm e}^{x_1x_2},$ and the Dirichlet boundary condition $h=-\frac{1}{\pi}\ln |\mathbf{x}|+{\rm e}^{x_1x_2}\cos x_1$. The analytic solution $u$ is given by
$u=-\frac{1}{\pi}\ln|\mathbf{x}|+{\rm e}^{x_1x_2}\cos x_1$.
\end{example}

\begin{figure}[H]
\centering \setlength{\tabcolsep}{0pt}
\begin{tabular}{ccc}
\includegraphics[height=4.2cm]{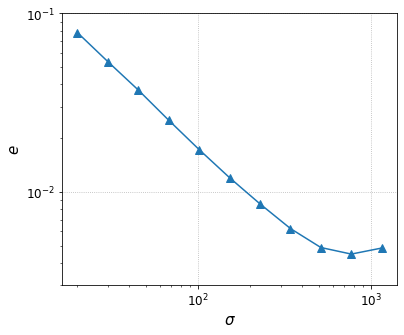}&
\includegraphics[height=4.2cm]{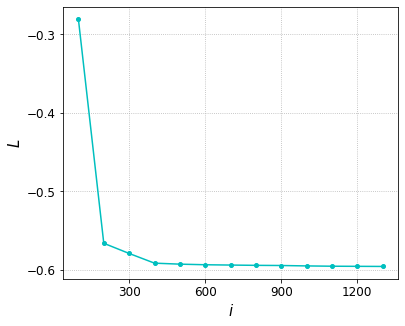}&
\includegraphics[height=4.2cm]{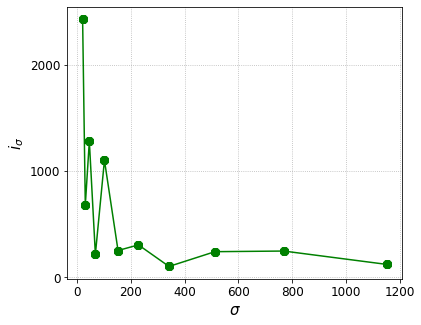}\\
(a) $e$ vs $\sigma$ & (b) $L$ vs $i$ & (c) $i_\sigma$ vs $\sigma$
\end{tabular}
\caption{\label{fig:Diri2d-conv} The dynamics of the training process: {\rm(a)} the error $e$ versus the penalty parameter $\sigma$, {\rm(b)} the decay of the loss function $L$ over iteration index $i$ with $\sigma = 1153.3$ and {\rm(c)} the required iteration number $i_\sigma$ for each $\sigma_k$-problem during the path following strategy.}
\end{figure}

\begin{figure}[H]
\centering \setlength{\tabcolsep}{0pt}
\begin{tabular}{cccc}
\includegraphics[height=2.8cm]{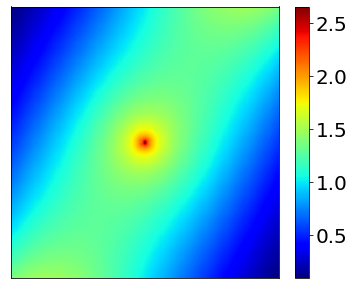}& \includegraphics[height=2.8cm]{true.png}& \includegraphics[height=2.8cm]{true.png}& \includegraphics[height=2.8cm]{true.png}\\
\includegraphics[height=2.8cm]{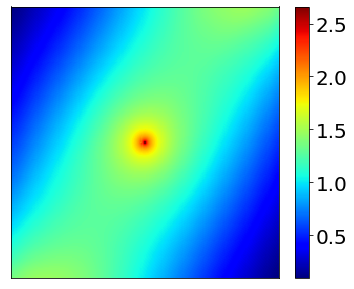} & \includegraphics[height=2.8cm]{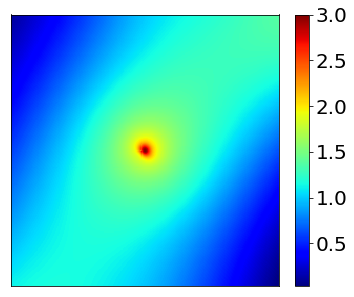}& \includegraphics[height=2.8cm]{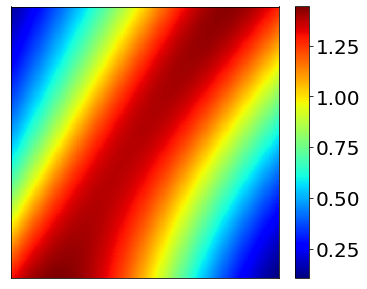} & \includegraphics[height=2.8cm]{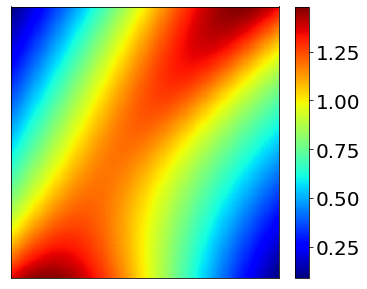}\\
\includegraphics[height=2.8cm]{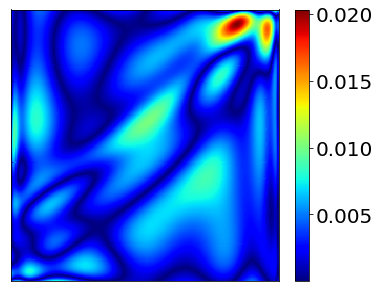} &
\includegraphics[height=2.8cm]{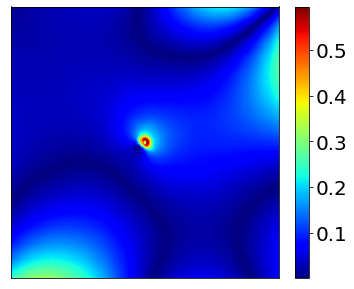} & \includegraphics[height=2.8cm]{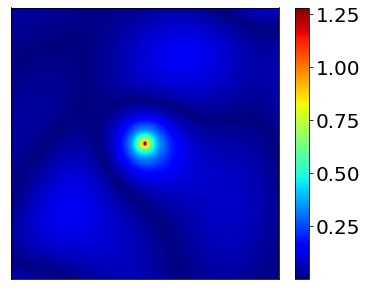} & \includegraphics[height=2.8cm]{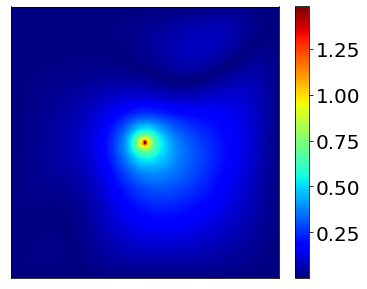}\\
SSDRM & SAPINN & DRM & WAN
\end{tabular}
\caption{The numerical approximations of Example \ref{exam:2d-1ps-Dirichlet} by the proposed SSDRM, SAPINN, DRM  and WAN. From top to bottom: analytic solution, neural network approximation and absolute error. \label{fig:exam:2d-1ps-Dirichlet}}
\end{figure}

In the splitting, the regular part $v$ is given by $v={\rm e}^{x_1x_2}\cos x_1$ and has to be learned using NNs.
The used NN has 3 hidden layers and each layer
has 20 neurons. The second condition is satisfied for $\sigma_{11}=1153.3$.
First we examine the convergence of the path-following strategy. We report the variation of the error $e$ with respect to the penalty parameter $\sigma$ in Fig. \ref{fig:Diri2d-conv}(a). It is observed that as the parameter $\sigma$ gets larger, the error $e$ drops rapidly during the first
few loops (i.e., $\sigma<400$), after which the decay becomes much slower. The decay plot roughly obeys the exponential law, which is consistent with the theoretical analysis in Section \ref{sec:theory}. Indeed, Theorem \ref{thm:error} indicates that as the value of $\sigma$ increases, the (mean squared) error due to penalization decreases to zero at a rate $O(\sigma^{-2})$, but the approximation error and the statistical error increase linearly with $\sigma$, which is also clearly observed in Fig. \ref{fig:Diri2d-conv}. Thus there is an optimal $\sigma$ which minimizes the total error. This plot shows also that one can terminate
the iteration earlier without compromising much the accuracy of the NN approximation. The empirical loss $\widehat {L}_\sigma$ decreases throughout the path-following procedure, cf. Fig. \ref{fig:Diri2d-conv}(b), exhibiting a steady convergence behavior of the optimizer. Note that the empirical loss for (SS)DRM is not necessarily positive (but it is bounded below by a constant).  Further, the number of iterations required for each $\sigma_k$-problem ($k\geq2$) is indeed much smaller than that for the $\sigma_1$-problem, confirming the effectiveness of the path-following strategy for warm-starting the optimization for $\sigma_k$-problems.

The numerical approximations by SSDRM and three existing NN based PDE solvers, i.e., SAPINN, DRM and WAN, are shown in Fig. \ref{fig:exam:2d-1ps-Dirichlet}, where the employed NN networks, the number of sampling points taken in the domain $\Omega$ and on the boundary $\partial\Omega$ and the
optimizer are identical for all the methods in order to ensure a fair comparison. The pointwise error of the SSDRM approximation is very small and the accuracy around the singularity at the origin is excellent, indicating a highly accurate approximation.
The results show that a direct application of DRM
fails to yield satisfactory results: the error is very large everywhere. This is not surprising since the corresponding continuous loss function $L_\sigma(v)$ is actually ill defined in the space $H^1(\Omega)$, as pointed out in Section \ref{ssec:existing}. Similar to DRM, WAN also does not perform well since it suffers from the exactly same issue.
In contrast, SAPINN can indeed yield much more accurate approximations than DRM and WAN, but the error around the singularity dominates the overall error. Intuitively, this concurs with the worse approximation property of NNs to singular functions \cite{DeVoreHanin:2021}. Table \ref{table:2} shows the quantitative comparison of these methods for
Example \ref{exam:2d-1ps-Dirichlet} in terms of computing time (in second) and error $e$. Note that the computing time should be interpreted indicatively only due to its sensitivity with respect to the setting and implementation etc. The results fully confirm the qualitative results in Fig. \ref{fig:exam:2d-1ps-Dirichlet}: DRM and WAN perform very poorly
on the example, and SAPINN does yield an acceptable approximation to the true solution $u$, but it is less accurate and computationally more expensive than SSDRM. It is worth mentioning that the accuracy of neural PDE solvers tends to stagnate at a level of $10^{-2}\sim 10^{-3}$, as for both SSDRM and SAPINN, but not much lower. This behavior has been widely observed for neural PDE solvers \cite{RAISSI2019686,yu2018deep,Zang:2020,CuomoSchiano:2022}, and it differs markedly from more traditional solvers, e.g., finite element method, for which one can make the error arbitrarily small by taking a sufficiently refined mesh. Also the computing time is relatively long, given the two-dimensional nature of the problem, indicating the imperative need for an accelerated training process, e.g., pretraining or meta-learning, in order for the method to be competitive with more traditional methods in the low-dimensional setting. We refer interested readers to \cite{PengKuhl:2021,RezaeiReese:2022} for detailed discussions on the pros and cons of neural PDE solvers and more classical approaches.

\begin{table}[H]
  \centering
  \begin{threeparttable}  \caption{The comparison between DRM, WAN, SAPINN, and the proposed SSDRM for Example \ref{exam:2d-1ps-Dirichlet}. \label{table:2}}
  \begin{tabular}{ccc}
    \toprule
    method & time (in sec) & $e$\\
    \midrule
    SSDRM& {252} & 4.85e-3\\
    SAPINN & {3433} & 7.82e-2\\
    DRM & {3} & 1.02e-1\\
    WAN & {32} & 1.38e-1\\
    \bottomrule
    \end{tabular}
    \end{threeparttable}
\end{table}

The next example is about the Poisson problem with a Neumann boundary condition.
\begin{example}\label{exam:2d-1ps-Neumann}
The domain $\Omega$ is the square
$\Omega=(-1,1)^2$, diffusion coefficient $\kappa\equiv1$, a point source at $\mathbf{x}_0=(0,0)$ with strength $f\equiv1$, source $g=\pi^2(\cos(\pi x_1)+\cos(\pi x_2))$, and
the Neumann boundary condition $h=\frac{\partial \Phi}{\partial n}$. The analytic solution $u$ is given by $u=\frac1{2\pi}\ln |\mathbf{x}| + \cos(\pi x_1)+\cos(\pi x_2)$.
\end{example}

\begin{figure}[H]
\centering\setlength{\tabcolsep}{0pt}
\begin{tabular}{ccc}
\includegraphics[height=4cm]{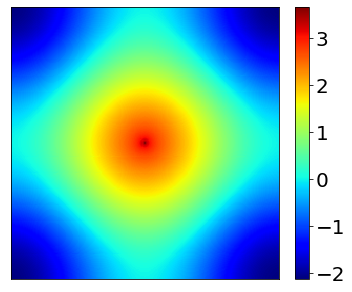} & \includegraphics[height=4cm]{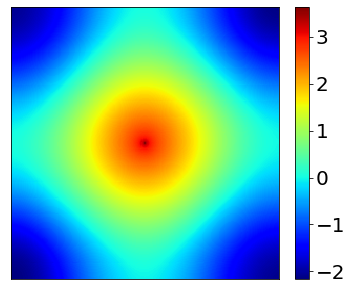} & \includegraphics[height=4cm]{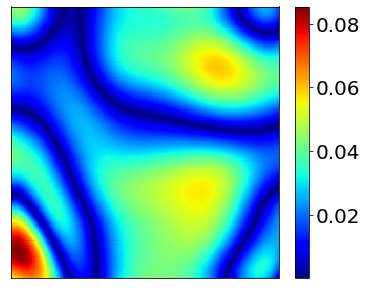}\\
 (a) exact & (b) SSDRM  & (c) error
\end{tabular}
\caption{\label{fig:exam:2d-1ps-Neumann}The numerical approximations of Example \ref{exam:2d-1ps-Neumann}, with an error $e=\text{3.18e-2}$.}
\end{figure}

In the splitting, the regular part $v(\mathbf{x})=\cos(\pi x_1)+\cos (\pi x_2)$ solves the standard Poisson problem with a source $F(\mathbf{x})=\pi^2(\cos(\pi x_1)+\cos(\pi x_2))$ and a zero Neumann boundary condition.
The value of the solution $v$ at a point $x^*$ on the boundary $\partial\Omega$ is specified to ensure the uniqueness of the exact solution. We employ an NN that has 3 layers and each layer has 8 neurons.
The obtained NN approximation after about 2500
iterations is shown in Fig. \ref{fig:exam:2d-1ps-Neumann} with an error $e=\textrm{3.18e-2}$. The observations from Example \ref{exam:2d-1ps-Dirichlet} are fully confirmed.

\begin{figure}[H]
\centering\setlength{\tabcolsep}{0pt}
\begin{tabular}{ccc}
 \includegraphics[height=4cm]{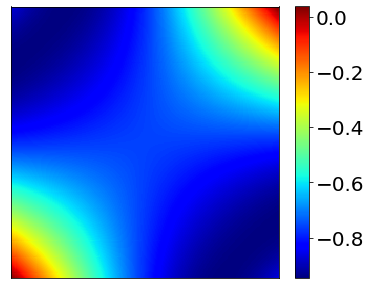} & \includegraphics[height=4cm]{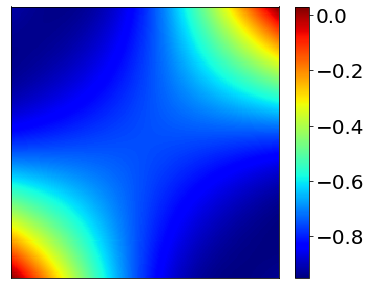} & \includegraphics[height=4cm]{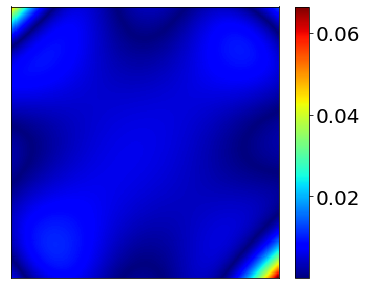}\\
 \includegraphics[height=4cm]{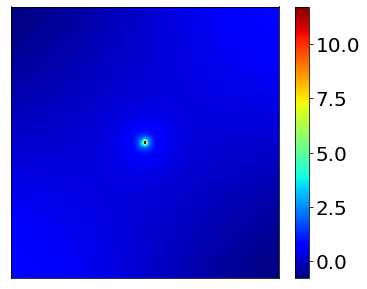} & \includegraphics[height=4cm]{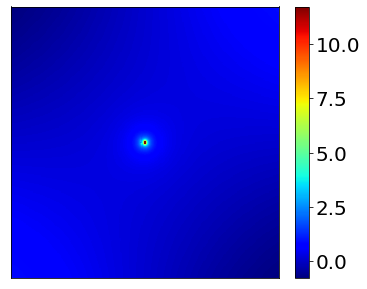} & \includegraphics[height=4cm]{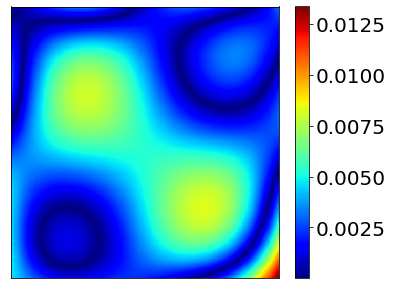}\\
 (a) exact & (b) SSDRM  & (c) error
\end{tabular}
\caption{\label{fig:m-d-point} The neural network approximation for Example \ref{exam:3d-1-sing}, with two slices at $x_3=-1$ (top) and $x_3=0$ (bottom).}
\end{figure}

Now we give a 3D problem with a single point source.
\begin{example} \label{exam:3d-1-sing}
The domain $\Omega$ is the cube $\Omega=(-1,1)^3$, diffusion coefficient $\kappa=|\mathbf{x}|^2+1$, a point source at the origin $\mathbf{x}_0=(0,0,0)$ with strength $f=|\mathbf{x}|^2+1,$ source $g=(|\mathbf{x}|^2+1)(|x_1^2+x_2^2+1)\sin(x_1x_2+x_3)-(4x_1x_2+2x_3)\cos(x_1x_2+x_3)+\frac{1}{2\pi|\mathbf{x}|}
$, and the Dirichlet boundary condition  $h=\frac{1}{4\pi|\mathbf{x}|}+\sin(x_1x_2+x_3)$.
The exact solution $u$ is given by $u=\frac{1}{4\pi|\mathbf{x}|}+\sin(x_1x_2+x_3).$
\end{example}
	
In the splitting, the regular part $v$ is given by $ v=\sin(x_1x_2+x_3) $. The used NN $v_\theta$ has 3 layers and each layer has 5 neurons, and the error $e$ of the SSDRM approximation is  $e=\text{8.14e-3}$. The numerical approximation is shown in Fig. \ref{fig:m-d-point},
with the two slices at $x_3=-1$ and $x_3=0$. The former illustrates the accuracy of boundary fitting, whereas
the latter shows that the approximation is of high quality near the singularity. The plots indicate that the pointwise error near the boundary $\partial\Omega$ tends to be slightly larger than that in the interior of the domain $\Omega$.

\begin{figure}[H]
\centering\setlength{\tabcolsep}{0pt}
\begin{tabular}{cc}
    \includegraphics[height=5.5cm]{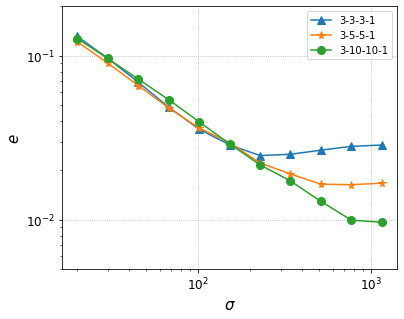} &
    \includegraphics[height=5.5cm]{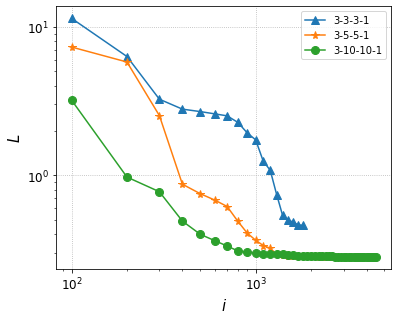} \\
    (a) $e$ vs $\sigma$ & (b) $L$ vs $i$
\end{tabular}
 \caption{\label{fig:m-d-point-penalty} The convergence of the path-following strategy for Example \ref{exam:3d-1-sing} with three different NNs: {\rm(a)} the error $e$ versus $\sigma$ and {\rm(b)} the empirical loss $\widehat{L}_\sigma$ versus the iteration index $i$ for $\sigma$ = 1153.3.}
\end{figure}

\begin{table}[H]
\begin{center}
\begin{threeparttable}
    \caption{The variation of the training time (in second) and error with respect to sampling points in the domain $N_r$ and on the boundary $N_{b}$.}\label{table:1}
\begin{tabular}  {cccccc}
\toprule
$N_{b}$ & $N_{r}$ & layer & neurons & time & $e$\\
\midrule
400 & 10000 & 3 & 3 & 45 & 2.86e-2\\
400 & 10000 & 3 & 5 & 80 & 1.70e-2\\
400 & 10000 & 3 & 10 & 130 & 9.66e-3\\
\midrule
400 & 10000 & 5 & 3 & 51 & 2.25e-2\\
400 & 10000 & 5 & 5 & 104 & 8.63e-3\\
400 & 10000 & 5 & 10 & 327 & 1.79e-2\\
\midrule
600 & 10000 & 3 & 3 & 47 & 3.45e-2\\
600 & 10000 & 3 & 5 & 83 & 8.14e-3\\
600 & 10000 & 3 & 10 & 132 & 9.60e-3\\
\midrule
400 & 15000 & 3 & 3 & 58 & 5.95e-2\\
400 & 15000 & 3 & 5 & 86 & 1.10e-2\\
400 & 15000 & 3 & 10 & 138 & 1.63e-2\\
\bottomrule
\end{tabular}
\end{threeparttable}
\end{center}
\end{table}

The variation of the error $e$ of the NN approximation as the penalty parameter $\sigma$ increases is shown in Fig. \ref{fig:m-d-point-penalty}.
In fact, as the penalty parameter $\sigma$ increases, the statistical and approximation
errors of the boundary term increase, whereas the error due to the penalization decreases. Thus, it is necessary to find a suitable value for $\sigma$ so as to achieve the optimal error of the NN approximation to the PDE \eqref{problem}. Fig. \ref{fig:m-d-point-penalty} indicates that one may choose $\sigma=1200$ to realize a good balance between difference sources of errors. Throughout, the value of the empirical loss $\widehat{L}_\sigma$ decays steadily to a limit. Also, we  have tested different numbers of sampling points and different architectures (layers, width) of the NN. The numerical results (i.e., the error $e$ and  the training time in seconds) are shown in Table \ref{table:1}. The results indicate that in order to achieve the best approximation with an optimal complexity, one needs a suitable balance between these different parameters. Increasing the size of the NN or the number of sampling points alone does not necessarily lead to a smaller error for the approximation. This is also observed earlier in \cite{JinLiLu:2022ip}. Thus, in the experiments, we have chosen the NN of the setting 3-5-5-1 (with $N_b=600$ and $N_r=10000$) to examine the training dynamics of the optimization method. Fig. \ref{fig:m-d-point-penalty} also indicates one peculiar but very common behavior of PDE solvers based on NNs: the accuracy of the NN approximation does not exhibit a steady decrease of the error as the NN width and depth increases, unlike the more traditional methods, e.g., finite element approximation, for which the error does decrease steadily to zero as the mesh size tends to zero. In sharp contrast, it appears that none of existing neural PDE solvers can exhibit a steady and consistent convergence rate in practical computation. The figure indicates that the iteration will stop at different functional values (with same penalty parameter $\sigma$) for the three NNs. For the 3-3-3-1 NN, the error is larger than the other two NNs, and the loss value cannot stabilize at a lower level. This is generally attributed to the optimization error associated with finding a global minimizer to the empirical loss $\widehat{L}_\sigma(v_\theta)$: instead the (stochastic) optimizer generally may find an approximate minimizer, due to the highly nonconvex landscape of the empirical loss function $\widehat{L}_\sigma$. The optimization error appears to be dominating for many DNN solvers \cite{JinLiLu:2022ip}, and its precise mathematical characterization is still largely open. We refer interested readers to \cite{Cyr:2020,AdcockDexter:2021} for detailed discussions on the persistent optimization error.

\begin{figure}[H]
\centering\setlength{\tabcolsep}{0pt}
\begin{tabular}{cc}
\includegraphics[width=7.5cm]  {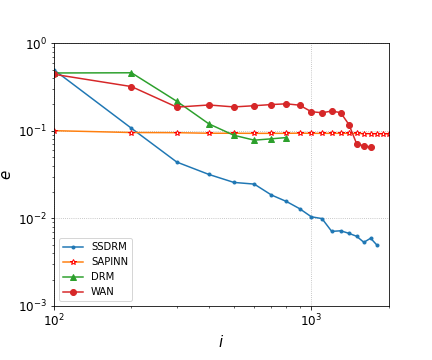} &
\includegraphics[width=7.5cm]  {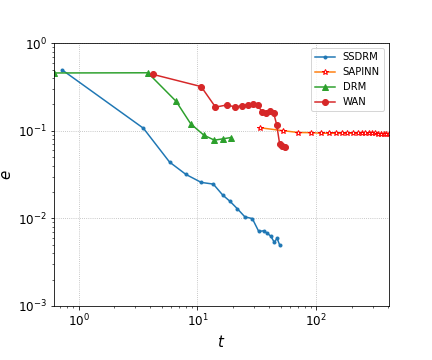}\\
(a) $e$ vs $i$ & (b) $e$ vs $t$
\end{tabular}
\caption{The evolution of the error versus the iteration index $i$ and the time $t$, respectively, for Example \ref{exam:3d-1-sing}, with SSDRM, SAPINN, DRM, and WAN. \label{fig:comp-3d}}
\end{figure}

Now we compare the proposed SSDRM with three existing NN based PDE solvers. See Fig. \ref{fig:comp-3d} for the error $e$ versus the iteration index $i$ and the computing time $t$ (in second). The proposed SSDRM significantly outperforms all existing approaches, and SAPINN beats both DRM and WAN by a large margin, since the latter two fail spectacularly to yield reasonable approximations. This is attributed to the presence of strong solution singularity, as discussed in Section \ref{ssec:existing}. These observations again fully confirm that from Example \ref{exam:2d-1ps-Dirichlet}.

The next example on a line source is taken from \cite[Example 5.2]{gjerde2019splitting}.
\begin{example}\label{exam:linesource}
The domain $\Omega$ is the unit cube $\Omega=(0,1)^3$, diffusion coefficient $\kappa\equiv 1$, a line source supported on the line segment connecting
$\mathbf{x}_1=(0,0,0.2)$  and $\mathbf{x}_2=(0,0,0.8)$ with strength $ f=x_3$, and $g=0 $, with a Dirichlet boundary condition $h=\frac{x_3}{4\pi}\ln\frac{|\mathbf{x}-\mathbf{x}_2|-(x_3-0.8)}{|\mathbf{x}-\mathbf{x}_1|-(x_3-0.2)}-(|\mathbf{x}-\mathbf{x}_1|-|\mathbf{x}-\mathbf{x}_2|)$. The exact solution $u$ is given by
$u=\frac{x_3}{4\pi}\ln\frac{|\mathbf{x}-\mathbf{x}_2|-(x_3-0.8)}{|\mathbf{x}-\mathbf{x}_1|-(x_3-0.2)}-(|\mathbf{x}-\mathbf{x}_1|-|\mathbf{x}-\mathbf{x}_2|).$
\end{example}

\begin{figure}[H]
   \centering
   \setlength{\tabcolsep}{0pt}
   \begin{tabular}{ccc}	\includegraphics[height=4cm]  {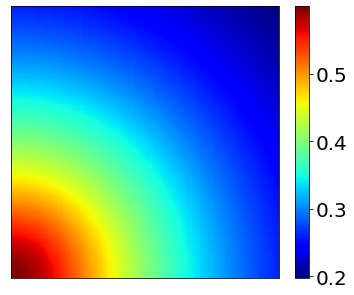} &
   \includegraphics[height=4cm]  {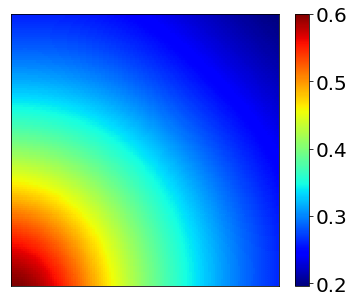} &		\includegraphics[height=4cm]  {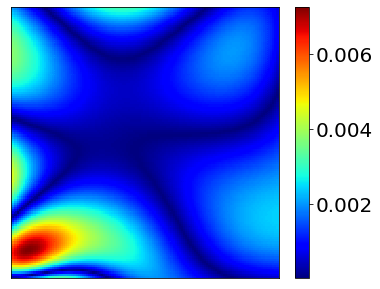}\\
   \includegraphics[height=4cm]  {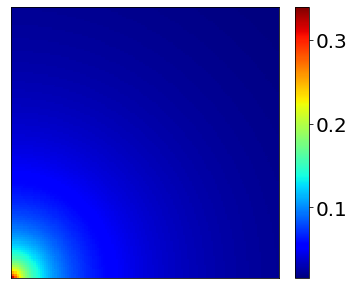} &
   \includegraphics[height=4cm]  {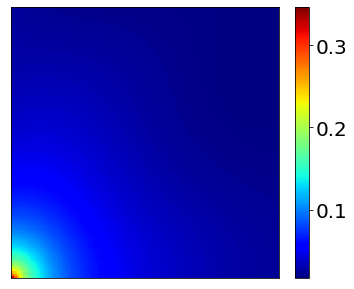} &
   \includegraphics[height=4cm]  {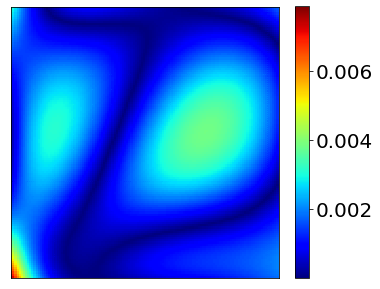}\\
   (a) exact & (b) SSDRM & (c) error
   \end{tabular}
\caption{\label{fig:linesource} The numerical
approximation for Example \ref{exam:linesource}, slice at $x_3=0$ (top) and $x_3=1/2$ (bottom).}
\end{figure}

\begin{figure}[H]
\centering\setlength{\tabcolsep}{0pt}
\begin{tabular}{cc}
\includegraphics[height=4.5cm]  {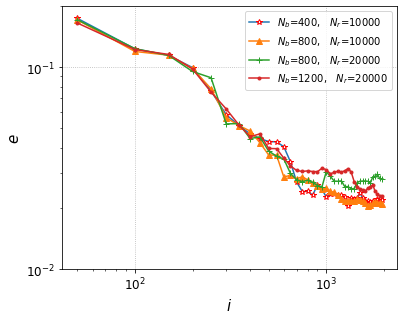}   &  \includegraphics[height=4.5cm]  {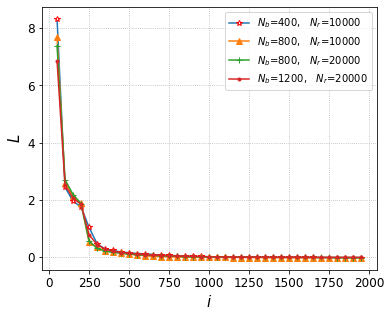}\\
\multicolumn{2}{c}{(a) architecture 3-10-10-1}\\
\includegraphics[height=4.5cm]  {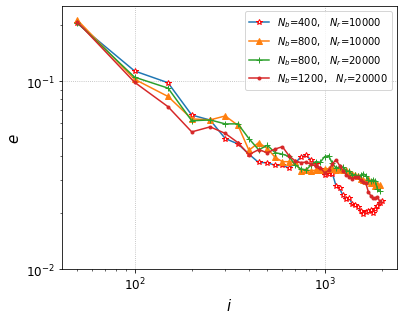}   &  \includegraphics[height=4.5cm]  {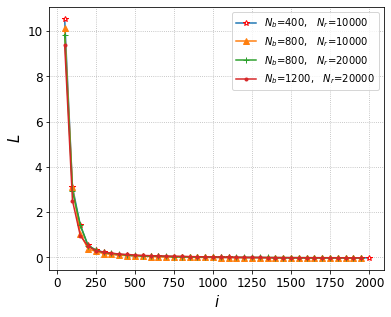}\\
\multicolumn{2}{c}{(b) architecture 3-10-10-10-1}\\
\end{tabular}
\caption{The evolution of the error $e$ and empirical loss $\widehat{L}_\sigma$ during the SSDRM training of Example \ref{exam:linesource} with different numbers $($i.e. $N_r$ and $N_b$$)$ of samples in the domain $\Omega$ and on the boundary $\partial\Omega$ for two neural network architectures.\label{fig:exam3d-line-conv}}
\end{figure}

In this example, we have to construct a function $\Phi_L$ whose Laplacian is a line source $\delta_{\Lambda}$ (supported on a line segment $\Lambda$). See the construction of $\Phi_L$ in Section \ref{ssec:split}.
Note that the singularities (lines) are actually placed on the boundary $\partial\Omega$ of the domain $\Omega$.
Then we split the solution into a singular part and a regular part $v$ and approximate $v$ by an NN $v_\theta$. In the experiment, we have employed an NN $v_\theta$ with 4 layers, each layer having 10 neurons. The stopping criterion is met for the penalty factor $\sigma=8000$, and the corresponding error is $e=\text{2.03e-2}$.  Two slices at $x_3=0$ and $x_3=1/2$ of the NN approximation are shown in Fig. \ref{fig:linesource}.
The slice at $x_3=0$ lies on the boundary $\partial\Omega$ of the domain $\Omega$. The
error on the boundary $\partial\Omega$ is
larger than that in the interior of the domain $\Omega$, but overall the result is still quite satisfactory. From the slice at $x_3=1/2$, one can clearly
identify the location of the line source (origin), and SSDRM can also yield good results around
the singular point.

\begin{figure}[H]
\centering
\setlength{\tabcolsep}{0pt}
\begin{tabular}{ccc}
\includegraphics[height=3.8cm]  {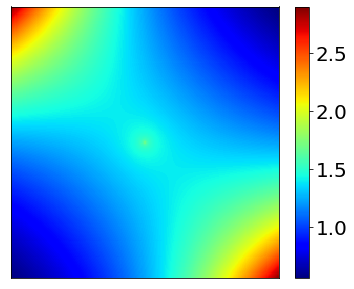} &
\includegraphics[height=3.8cm]  {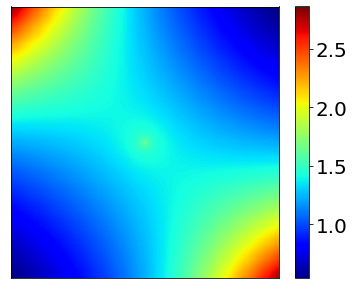} &
\includegraphics[height=3.8cm]  {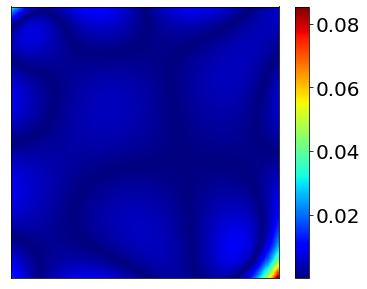}\\
\includegraphics[height=3.8cm]  {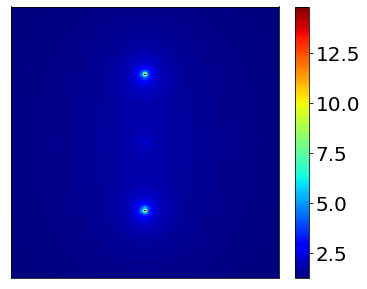} &
\includegraphics[height=3.8cm]  {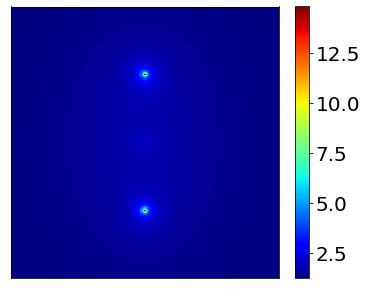} &
\includegraphics[height=3.8cm]  {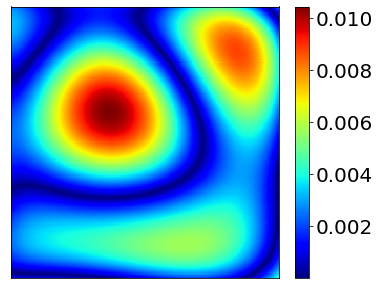}\\
(a) exact & (b) SSDRM & (c) error
\end{tabular}
\caption{\label{Figure 9}The neural network approximation for Example \ref{exam:multi-sing}, slices at $x_3=-1$ (top) and $x_3=0$ (bottom).}
\end{figure}

To examine the influence of the numbers of sample points (i.e., $N_r$ and $N_b$) and the NN architecture on the convergence behavior of the training process (in the sense of optimization), we have experimented with different settings. The results are presented in Fig. \ref{fig:exam3d-line-conv}, where the penalty parameter $\sigma$ is fixed at $\sigma=8000$.
The training of SSDRM converges robustly for all settings: the loss value decreases very steadily, and for all the settings the approximations exhibit very much comparable accuracy and achieve similar loss values. Although not presented, it is noted that the training may suffer from ``over-fitting'', where the error $e$ starts to increase after some initial iterations, when the penalty parameter $\sigma$ is much smaller, e.g., $\sigma=2000$, $N_{b}=400$ and $N_r = 10000$ for the architecture 3-10-10-10-1.
This agrees with previous discussions that one needs to choose the penalty parameter $\sigma$ properly in order to achieve the best possible approximation.

The proposed SSDRM extends straightforwardly to the more complex case involving more than one
singularity, or a combination of point and line sources. The next three-dimensional example
contains a line source and two point singularities.
\begin{example}\label{exam:multi-sing}
The domain $\Omega$ is $\Omega=(-1,1)^3$, diffusion coefficient $ \kappa={\rm e}^{x_3} $, with a line segment source $\{(0,0,x_3): -1\le x_3\leq 1\}$ with two endpoints $\mathbf{x}_1=(0,0,-1)$, $\mathbf{x}_2=(0,0,1)$, and two point sources concentrated at $\mathbf{x}_3=(0,\frac{1}{2},0)$ and $
\mathbf{x}_4=(0,-\frac{1}{2},0)$ and  $f={\rm e}^{x_3}$ and $g=
-(x_1x_2+x_1^2x_2^2+x_2^2x_3^2+x_3^2x_1^2){\rm e}^{x_1x_2x_3+x_3}+\frac{x_3{\rm e}^{x_3}}{4\pi|\mathbf{x}-\mathbf{x}_3|^3}+\frac{x_3{\rm e}^{x_3}}{4\pi|\mathbf{x}-\mathbf{x}_4|^3}+\frac{{\rm e}^{x_3}}{4\pi}(\frac{1}{|\mathbf{x}-\mathbf{x}_2|}-\frac{1}{|\mathbf{x}-\mathbf{x}_1}|),$ and the Dirichlet boundary condition
$h=\frac{1}{4\pi}\ln\frac{|\mathbf{x-x_2}|-(x_3-1)}{|\mathbf{x-x_1}|-(x_3+1)}+\frac{1}{4\pi}\frac{1}{|\mathbf{x}-\mathbf{x}_3|}+
\frac{1}{4\pi}\frac{1}{|\mathbf{x}-\mathbf{x}_4|}+{\rm e}^{x_1x_2x_3}.$
The analytic solution $u$ is given by
$u=\frac{1}{4\pi}\ln\frac{|\mathbf{x-x_2}|-(x_3-1)}{|\mathbf{x-x_1}|-(x_3+1)}+\frac{1}{4\pi}\frac{1}{|\mathbf{x}-\mathbf{x}_3|}+
\frac{1}{4\pi}\frac{1}{|\mathbf{x}-\mathbf{x}_4|}+{\rm e}^{x_1x_2x_3}.$
\end{example}

The NN architecture for this example is taken to be 3-6-6-1. The training terminates at $e=\text{4.37e-3}$ and $\sigma=1730$. The NN approximation is shown in Fig. \ref{Figure 9}, with the two slices at $x_3=-1$ and $x_3=0$, respectively. The results confirm previous observations: the singular part is exactly resolved, and the regular part $v$ is well approximated so that SSDRM can produce an excellent approximation to the exact solution. This example shows clearly the flexibility of the proposed approach for elliptic problems with multiple singularities.

The next example further shows the flexibility of the proposed SSDRM
for resolving singular sources living on high-dimensional
subspaces, as mentioned in Section \ref{ssec:split}.
\begin{example}\label{exam:5d}
This example is 5D Poisson equation, with the domain $\Omega=(-1,1)^5$,
surface sources at $\mathbf{x}_1'=(\frac12,\frac12),$ $\mathbf{x}_2'=(0,\frac12)$, $\mathbf{x}_3'=(-\frac12,\frac12)$,
$\mathbf{x}_4'=(\frac12,0)$, $\mathbf{x}_5'=(0,0)$, $\mathbf{x}_6'=(-\frac12,0)$, $\mathbf{x}_7'=(\frac12,-\frac12)$, $\mathbf{x}_8'=(0,-\frac12)$,
and $\mathbf{x}_9'=(-\frac12,-\frac12)$ {\rm(}with $\mathbf{x}'=(x_1,x_2)${\rm)}, $f\equiv1$ and $g=-10$. The Dirichlet boundary
condition $h=\sum_{i=1}^9-\frac{1}{2\pi}\ln |\mathbf{x}'-\mathbf{x}_i'| + |\mathbf{x}|^2$. The analytic solution $u$ is given
by $u=\sum_{i=1}^9-\frac{1}{2\pi}\ln |\mathbf{x}'-\mathbf{x}_i'| + |\mathbf{x}|^2$.
\end{example}

\begin{figure}[H]
\centering
\setlength{\tabcolsep}{0pt}
\begin{tabular}{ccc}
\includegraphics[height=4cm]  {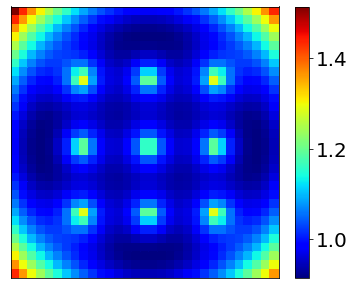} & \includegraphics[height=4cm]  {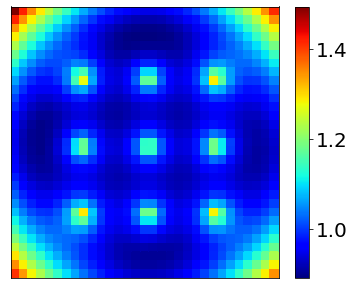} & \includegraphics[height=4cm]  {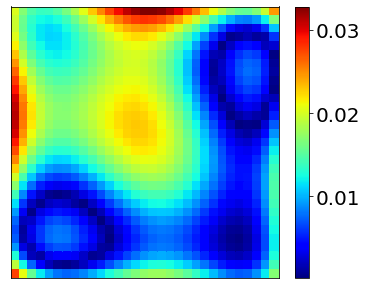}\\
(a) exact & (b) SSDRM & (c) error
\end{tabular}
\caption{\label{Figure 16}
The neural network approximation for Example \ref{exam:5d}, slice at $x_3=x_4=x_5=0$.}
\end{figure}

The singularity is understood to be supported on  three-dimensional affine subspace through the bounded region $\Omega$, similar to Gjerde et al \cite{gjerde2019splitting}. Then the corresponding singularity function capturing the source term $\delta(x_1',x_2')$ in 5-D space is given by $\Phi(\mathbf{x}')=-\frac{1}{2\pi}\ln |\mathbf{x}'-\mathbf{x}_i'|$. The NN for approximating the regular part $v$ has 2 layers and each hidden layer has 10 neurons. The loop terminates at the penalty factor $\sigma=1922$, and the error $e$ of the NN approximation is $e=\text{8.14e-3}$ in the domain. One slice of the approximation at $x_3=x_4=x_5=0$ is shown in Fig. \ref{Figure 16}. The numerical results show again the singularities are accurately resolved with a small pointwise error.

\section{Error analysis}\label{sec:theory}
Now we discuss  the error analysis for SSDRM. Due to the use of singularity splitting, the issue boils down to the analysis of the standard DRM for smooth solutions, which has been carried out independently in several recent works
under different problem settings and architectures \cite{LuLuWang:2021,DuanJiao:2021cicp,duan2021analysis,HongSiegelXu:2021,Muller:2021Error,LuChenLu:2021}. Below we summarize relevant theoretical results for the Dirichlet case, to provide theoretical underpinnings of the SSDRM and also shed valuable insights into the influence of various algorithmic parameters (e.g., penalty parameter $\sigma$ and NN width and depth). The Neumann case is simpler (since there is no error associated with the penalization) and can be treated similarly. Throughout the discussion, the domain $\Omega$ is assumed to be smooth.

\subsection{Error estimate of the penalty method}
First, we discuss the error incurred by the penalty method, i.e., how the error decreases as the penalty
parameter $\sigma\to \infty$.  We denote the minimizer of $L_{\sigma}(v)$ by $v_\sigma$ and
the solution of problem \eqref{prob} by $v^*$. The next result shows that $v_\sigma$
actually can be viewed as a Robin approximation to the Dirichlet problem.
\begin{proposition}
Let $v_\sigma\in H^1(\Omega)$ be the minimizer of $L_\sigma(v)$. Then $v_\sigma$ is the weak solution of the following Robin boundary value problem:
\begin{equation}\label{app prob}
\left\{\begin{aligned}
	-\nabla\cdot(\kappa\nabla v)&=F, &&\mbox{in }\Omega,\\
	v+\frac{\kappa}{\sigma}\frac{\partial v}{\partial n}&=\tilde h,&& \mbox{on }\partial\Omega.
\end{aligned}\right.
\end{equation}
\end{proposition}
\begin{proof}
Since $v_\sigma\in H^1(\Omega)$ is the minimizer of $L_\sigma(v)$, for every $ \phi\in H^1(\Omega) $ we have $ L_\sigma(v_\sigma)\leq L_\sigma(v_\sigma+\lambda \phi) $, $\lambda\in\mathbb{R}$. Let $ f(\lambda)=L_\sigma(v_\sigma+\lambda \phi) $.
Direct computation gives
\begin{equation*}
f'(\lambda)=\lambda\left((\kappa\nabla \phi,\nabla \phi)+\sigma(\phi,\phi)_{L^2(\partial\Omega)}\right)+(\kappa\nabla v_\sigma,\nabla \phi)-(F,\phi)+\sigma(v_\sigma-\tilde h,\phi)_{L^2(\partial\Omega)}.
\end{equation*}
The minimizing property of $v_\sigma$
implies $ f'(0)=0 $, and consequently
\begin{equation}\label{variation}
f'(0)=(\kappa\nabla v_\sigma,\nabla \phi)-(F,\phi)+\sigma(v_\sigma-\tilde h,\phi)_{L^2(\partial\Omega)}=0,\quad\forall \phi\in H^1(\Omega).
\end{equation}
This is the variational formulation of \eqref{app prob}. So $v_\sigma$ is the weak solution of problem \eqref{app prob}.
\end{proof}

The penalization is a type of singular perturbation. It has been studied under various situations \cite{CostabelDauge:1996,Auchmuty:2018}. To derive error bounds, we recall a regularity result for the Robin boundary value problem essentially due to Costabel and Dauge \cite{CostabelDauge:1996}.
\begin{lemma}\label{lemma:Costabel}
Let $\kappa\in C^2(\overline{\Omega})$. Fix $\epsilon>0$, and let $z^\epsilon$ be the solution to the Robin boundary value problem
\begin{equation}\label{eqn:z}
   \left\{\begin{aligned}
    - \nabla \cdot(\kappa\nabla z^\epsilon) & = 0,\quad \mbox{in }\Omega,\\
    \epsilon \kappa\frac{\partial z^\epsilon}{\partial n} + z^\epsilon & = g,\quad \mbox{on }\partial\Omega.
 \end{aligned}\right.
\end{equation}
Then for any $t\in[0,3/2]$ and $t \le s\le t+1$, there exists some $C>0$ independent of $\epsilon$ such that
\begin{equation*}
    \|z^\epsilon\|_{H^{1+s}(\Omega)}\leq C\epsilon^{t-s}\|g\|_{H^{\frac{1}{2}+t}(\partial\Omega)}.
\end{equation*}
\end{lemma}
\begin{proof}
The proof follows closely the argument of \cite[Lemma 2.1 and Corollary 2.2]{CostabelDauge:1996} for the Laplacian, but it is adapted to the case of a variable coefficient.
The proof utilizes the Dirichlet-to-Neumann map: for a given $g\in H^{\frac12}(\partial\Omega)$, let
$\Lambda_\kappa g = \kappa \frac{\partial w_g}{\partial n}\big|_{\partial\Omega},$
where $w_g \in H^1(\Omega)$ solves
\begin{equation}\label{eqn:wg}
   \left\{\begin{aligned}
    - \nabla \cdot(\kappa\nabla w_g) & = 0,\quad \mbox{in }\Omega,\\
    w_g & = g,\quad \mbox{on }\partial\Omega.
 \end{aligned}\right.
\end{equation}
Let $T:H^s(\Omega) \rightarrow H^{s-\frac12}(\partial\Omega)$ be the trace operator. Then $Tz^\epsilon$ can be represented as
\begin{equation}\label{eqn:zep}
 Tz^\epsilon=(\epsilon \Lambda_\kappa+I)^{-1}g.
\end{equation}
Then it suffices to show the (uniform in $\epsilon$) boundedness of the operator $(\epsilon\Lambda_\kappa + I)^{-1}$ in $H^s(\partial\Omega)$ with $s\ge0$.
Note that for any $g\in H^{\frac12}(\partial\Omega)$
\begin{equation*}
  ((\epsilon\Lambda_\kappa + I) g, g)_{L^2(\partial\Omega)} = \epsilon \|\sqrt{\kappa}\nabla w_g \|_{L^2(\Omega)}^2 + \| g \|_{L^2(\partial\Omega)}^2 \ge  \|g\|_{L^2(\partial\Omega)}^2.
\end{equation*}
Therefore the operator $\epsilon\Lambda_\kappa + I$ is
 positive definite in $L^2(\partial\Omega)$ and
$$\|(\epsilon\Lambda_\kappa + I)^{-1}\|_{L^2(\partial\Omega)\rightarrow L^2(\partial\Omega)} \le C,$$
where the constant $C$ is independent of $\epsilon$. To show the uniform bound of the operator in $H^s(\partial\Omega)$, let $(x',x_n)$ be local coordinate near a point
$x_0 \in \partial\Omega$ so that the boundary is given by
$x_n = 0$. The function $\lambda_\kappa (x',\xi)$ denotes the
symbol of $\Lambda_\kappa$ in this coordinate. Then the following asymptotic expansion holds (see e.g., \cite[Theorem 4.1]{Uhlmann:2009} and \cite[Theorem 0.1]{SylvesterUhlmann:1988})
$$\lambda_\kappa(x',\xi') =  \kappa(x',0)| \xi'|  + a_0(x',\xi') + r(x',\xi'),$$
where $\kappa(x',0)|\xi'|$ denotes a classical pseudo-differential operator of order $1$, $a_0(x',\xi')$ is homogeneous of degree $0$ in $\xi'$,
and the residue $r(x',\xi')$ is a classical symbol of order $-1$.
Hence, for any $s\ge0$
\begin{equation*}
\begin{aligned}
&\| (\epsilon\Lambda_\kappa + I)^{-1} g \|_{H^{s}(\partial\Omega)} \le
C  \| (\Lambda_\kappa + I)^s (\epsilon\Lambda_\kappa + I)^{-1} g \|_{L^2(\partial\Omega)} \\
=&
C  \| (\epsilon\Lambda_\kappa + I)^{-1} (\Lambda_\kappa + I)^s g \|_{L^2(\partial\Omega)}
\le C \| (\Lambda_\kappa + I)^s g \|_{L^2(\partial\Omega)}  \le C \|g \|_{H^s(\partial\Omega)},
\end{aligned}
\end{equation*}
since $\Lambda_\kappa$ is self-adjoint and semi-positive definite, and $(\epsilon\Lambda_\kappa+I)^{-1}$ commutes with $(\Lambda_\kappa+I)^s$.
Then by the defining identity \eqref{eqn:zep}, we obtain
\begin{equation*}
 \| Tz^\epsilon \|_{H^{s}(\partial\Omega)} \le   C \|g \|_{H^{s}(\partial\Omega)}.
\end{equation*}
Now the smoothing property of the Dirichlet problem implies
\begin{equation}\label{eqn:z-1}
 \|  z^\epsilon \|_{H^{1+t}( \Omega)} \le C\| Tz^\epsilon \|_{H^{\frac12+t}(\partial\Omega)} \le   C \|g \|_{H^{\frac12+t}(\partial\Omega)},\quad \forall 0 \le t\le \tfrac32.
\end{equation}
Meanwhile, the Robin problem \eqref{eqn:z} can be rewritten as
	\begin{equation}
		\left\{\begin{aligned}
			- \nabla \cdot(\kappa\nabla z^\epsilon) & = 0,\quad \mbox{in }\Omega,\\
			\kappa\frac{\partial z^\epsilon}{\partial n}+z^\epsilon & = \epsilon^{-1}(g - (1-\epsilon) z^\epsilon),\quad \mbox{on }\partial\Omega.
		\end{aligned}\right.
	\end{equation}
Then  the smoothing property for
the Robin boundary value problem leads to \cite[Theorem 3.10 (iii)]{ErnGuermond:2004}
\begin{equation}\label{eqn:z-2}
\begin{aligned}
&\quad\|z^\epsilon\|_{H^{\frac32+t}(\Omega)}\leq C\epsilon^{-1}\|g - (1-\epsilon) Tz^\epsilon\|_{H^{\frac12+t}(\partial\Omega)}\\
&\leq C\epsilon^{-1}(\|g\|_{H^{\frac12+t}(\partial\Omega)} + \| T z^\epsilon \|_{H^{\frac12+t}(\partial\Omega)}) \leq C\epsilon^{-1}||g||_{H^{\frac12+t}(\partial\Omega)},
\end{aligned}
\end{equation}
where the last step is due to the estimate \eqref{eqn:z-1}.
Then interpolating between \eqref{eqn:z-1} and \eqref{eqn:z-2} gives the desired result.
\end{proof}

Now we can state a bound on the error $v_\sigma-v^*$, arising from approximating the Dirichlet problem with a Robin one.
\begin{theorem}\label{thm:penalization}
Let $F\in H^{t}(\Omega)$, and $\tilde h\in H^{t+\frac32}(\partial\Omega)$, with $t\in [-1,0]$. Then for any $t-1\leq s\leq t$, the following error bound holds
\begin{equation*}
 \|v_\sigma - v^*\|_{H^{1+s}(\Omega)} \leq     C\sigma^{s-t-1} (\|F\|_{H^t(\Omega)}+\|\tilde h\|_{H^{t+\frac32}(\partial\Omega)}).
\end{equation*}
\end{theorem}
\begin{proof}
Note that $w_\sigma = v_\sigma-v^*$ satisfies
\begin{equation*}
   \left\{\begin{aligned}
    - \nabla \cdot(\kappa\nabla w_\sigma) & = 0,\quad \mbox{in }\Omega,\\
    \sigma^{-1} \kappa\frac{\partial w_\sigma}{\partial n} + w_\sigma & =- \sigma^{-1}\kappa\frac{\partial v^*}{\partial n},\quad \mbox{on }\partial\Omega.
 \end{aligned}\right.
\end{equation*}
Then by Lemma \ref{lemma:Costabel} and the trace theorem, there holds for $t-1\leq s\le t$
\begin{align*}
   & \|w_\sigma\|_{H^{1+s}(\Omega)} \leq C\sigma^{s-t}\Big\|\sigma^{-1}\kappa\frac{\partial v^*}{\partial n}\Big\|_{H^{\frac{1}{2}+t}(\partial\Omega)}\\
    =&C\sigma^{s-t-1}\Big\|\kappa\frac{\partial v^*}{\partial n}\Big\|_{H^{\frac{1}{2}+t}(\partial\Omega)}
    \leq C\sigma^{s-t-1}\| v^*\|_{H^{2+t}(\Omega)}.
\end{align*}
Meanwhile, by the standard elliptic regularity theory, we have
\begin{equation*}
    \| v^*\|_{H^{2+t}(\Omega)} \leq C\Big(\|F\|_{H^t(\Omega)} + \|\tilde h\|_{H^{\frac32+t}(\partial\Omega)}\Big).
\end{equation*}
Combining the preceding two estimates completes the proof the theorem.
\end{proof}
\begin{remark}
Theorem \ref{thm:penalization} indicates that under suitable regularity assumptions on the solution $v^*$, the error due to penalization decays {\rm(}sub{\rm)}linearly.
The error bound in Theorem \ref{thm:penalization} is known in various special forms. For example, in both \cite{duan2021analysis} and \cite{Muller:2021Error}, the following error bound was derived for the case of $\tilde h\equiv0$ and $F\in L^2(\Omega)$:
\begin{equation*}
    \|v^*-v_\sigma\|_{H^1(\Omega)}\leq C\sigma^{-1}\|F\|_{L^2(\Omega)}.
\end{equation*}
It was derived for a strongly coercive bilinear form in the space $H^1(\Omega)$ in \cite{duan2021analysis}, using an energy argument. In \cite[Theorem 9]{Muller:2021Error}, it was derived using the theory of the Steklov eigenvalue problem {\rm(}see \cite{Auchmuty:2018} and references therein{\rm)}. In contrast to existing results, Theorem \ref{thm:penalization} gives the bound for a broader range of regularity conditions on the problem data $F$ and $\tilde h$.

Under the \textit{a priori} assumption $v^*\in H^1(\Omega)$ only, it can be proved
\begin{equation*}
    \|v_\sigma-v^*\|_{L^2(\Omega)}\leq C(\kappa)(\|F\|_{(H^{\frac12}(\Omega))'} + \|\tilde h\|_{H^\frac12(\Omega)})\sigma^{-\frac12}.
\end{equation*}
See \cite[Lemma 17]{Muller:2021Error} for a detailed proof when $\tilde h\equiv0$, and the case $\tilde h\not\equiv0$ follows similarly. One can also show convergence in $H^1(\Omega)$ but without a rate, using a standard compactness argument from calculus of variation, i.e.,
\begin{equation*}
    \lim_{\sigma\to\infty}\|v^*-v_\sigma\|_{H^1(\Omega)}=0.
\end{equation*}
Further discussions on the convergence of the penalization under even weaker regularity assumption for the homogeneous Laplace equation {\rm(}i.e., $F\equiv0${\rm)}, including very weak solutions with Dirichlet boundary data $\tilde h \in H^s(\Omega)$ with $0<s<\frac12$, can be found in \cite[Section 7]{Auchmuty:2018}
\begin{equation*}
    \| v_\sigma-v^*\|_{L^2(\Omega)} \leq C\sigma^{-s}\|\tilde h\|_{H^s(\partial\Omega)};
\end{equation*}
See \cite[Theorem 7.1]{Auchmuty:2018} for the precise statement.
\end{remark}

\subsection{Generalization error analysis}
In practical computation,  we use an empirical loss $\widehat{L}_\sigma(v_\theta)$ instead of the continuous loss $L_\sigma(v_\theta)$, due to the use of the Monte Carlo method to approximate the integrals. Let $\widehat{v_\sigma}$ be a (global) minimizer
of the empirical loss $\widehat{L}_\sigma(v_\theta)$ over the set $\mathcal{A}$ of neural network functions (with a fixed architecture, e.g., depth, number of nonzero parameters, and maximum bound on the parameter vectors). Now we bound the error between
$\widehat{v_\sigma}$ and $v_{\sigma} $. This kind of analysis is commonly known as the generalization error in statistical learning theory \cite{AnthonyBartlett:1999,ShalevShwartzBenDavid:2014}. For the discussion below, we assume the modified source $F\in L^\infty(\Omega),\tilde h \in L^\infty(\partial\Omega)$, which holds for $\kappa,f\in W^{2,\infty}(\Omega)$ and $h\in L^\infty(\partial\Omega)$ and being locally constant in a neighborhood of the singularity $\delta$. This assumption is due to the limitation of the proof technique using Rademacher complexity. The next lemma gives an error decomposition into the quadrature error and approximation error.

\begin{lemma}\label{lem:decom}
Let $\mathcal{A}$ be a set of NN functions of fixed depth, width, and upper bound on the parameters. Let $\widehat{v}_\sigma$ be the minimizer of $\widehat{L}_\sigma(v)$ over $\mathcal{A}$. Then
\begin{equation}\label{errorest}
\|\nabla(\widehat{v_\sigma}-v_\sigma)\|_{L^2(\Omega)}^2+\sigma\|\widehat{v_\sigma}-v_\sigma\|_{L^2(\partial\Omega)}^2\leq C(\kappa)(\sup_{\widehat{v}\in\mathcal{A}}|L_\sigma(\widehat{v})-\widehat{L_\sigma}(\widehat{v})|+\inf_{\widehat{v}\in\mathcal{A}}\sigma \|\widehat{v}-v_\sigma\|_{H^1(\Omega)}^2).
\end{equation}
\end{lemma}
\begin{proof}
This result follows from direct computation. For every $\widehat{v}\in\mathcal{A}$, by the minimizing property $ \widehat{L}_\sigma(\widehat{v_\sigma})
\leq\widehat{L}_\sigma(\widehat{v}) $, we have
\begin{align*}
L_\sigma(\widehat{v_\sigma})-L_\sigma(v_\sigma)		 =&[L_\sigma(\widehat{v_\sigma})-\widehat{L}_\sigma(\widehat{v_\sigma})]+[\widehat{L}_\sigma(\widehat{v_\sigma})-\widehat{L}_\sigma(\widehat{v})]\\
 &+[\widehat{L}_\sigma(\widehat{v})-L_\sigma(\widehat{v})]+[L_\sigma(\widehat{v})-L_\sigma(v_\sigma)]\\
\leq&2\sup_{\widehat{v}\in\mathcal{A}}|L(\widehat{v})-\widehat{L}_\sigma(\widehat{v})|+L_\sigma(\widehat{v})-L_\sigma(v_\sigma).
\end{align*}
Taking the infimum over $\widehat{v}\in\mathcal{A}$ gives
\begin{equation}\label{errorde}
L_\sigma(\widehat{v_\sigma})-L_\sigma(v_\sigma)\leq2\sup_{\widehat{v}\in\mathcal{A}}|L_\sigma(\widehat{v})-\widehat{L_\sigma}(\widehat{v})|+\inf_{\widehat{v}\in\mathcal{A}}[L_\sigma(\widehat{v})-L_\sigma(v_\sigma)].
\end{equation}
For any $\widehat{v}\in\mathcal{A}$, and $w=\widehat v-v_\sigma$, we deduce
\begin{align*}
L_\sigma(\widehat{v})-L_\sigma(v_\sigma)
=&\tfrac{1}{2}(\kappa\nabla w,\nabla w)+(\kappa\nabla v_\sigma,\nabla w)+\tfrac{\sigma}{2}\|w\|_{L^2(\partial\Omega)}^2\nonumber\\
&+\sigma(v_\sigma-\tilde h,w)_{L^2(\partial\Omega)}-(F,w)\nonumber\\
=&\tfrac{1}{2}(\kappa \nabla w,\nabla w)+\tfrac{\sigma}{2}\|w\|_{L^2(\partial\Omega)}^2,
\end{align*}
where the last step is due to the weak formulation of $v_\sigma$, i.e.,
\begin{equation*}
    (\kappa\nabla v_\sigma,\nabla w) + \sigma(v_\sigma-\tilde h,w)_{L^2(\partial\Omega)}=(F,w).
\end{equation*}
Upon substituting the identity into \eqref{errorde} and applying the trace theorem, we obtain the desired  bound.
\end{proof}

The error decomposition in Lemma \ref{lem:decom} provides interesting insights into the total error of the NN approximation $\widehat{v_\theta}$: The two terms $ \sup_{\widehat{v}\in\mathcal{A}}|L_\sigma(\widehat{v})-\widehat{L_\sigma}(\widehat{v})|$ and  $\inf_{\widehat{v}\in\mathcal{A}}\sigma ||\widehat{v}-v_\sigma||_{H^1(\Omega)}$ respectively represent the statistical error, due to the use of the Monte Carlo method to approximate the integral and the approximation error, due to the use of NNs to approximate the function $v_\sigma$.

Note that in the analysis, we have assumed that a global optimizer $\widehat{v_\sigma}$ of the empirical loss $\widehat{L}_\sigma(v_\theta)$ can be numerically realized, which is unfortunately generally not the case in practice. Indeed, the landscape of the empirical loss $\widehat{L}_\sigma(v_\theta)$ is highly nonconvex (see \cite{Krishnapriyan:2021} for graphical illustrations of PINN), due to the nonlinearity of the activation function $\rho$, and there is no guarantee that an optimizer will find a global minimizer $\widehat{\theta}$. Instead only an approximate minimizer $\tilde \theta$ can be expected in practice. This leads to an additional source of error, known as optimization error in the literature. It is often observed to be dominating in many neural PDE solvers \cite{WangPerdikaris:2022jcp,Krishnapriyan:2021,JinLiLu:2022ip}. Up to now, it is still completely open to rigorously analyze the optimization error, and thus is often ignored from the analysis.

Below we bound the approximation error and statistical error. This kind of analysis has been pursued in several recent works \cite{LuLuWang:2021,DuanJiao:2021cicp,JiaoLai:2021elliptic,Muller:2021Error,LuChenLu:2021} under different assumptions on the activation function $\rho$ and NN architecture. Under the assumption $F\in L^\infty(\Omega)$, the work \cite{JiaoLai:2021elliptic} provides relevant estimates on the quadrature error by means of the Monte Carlo method using Rademacher complexity, and the approximation error of the penalized solution with NNs using the approximation theory from \cite{I2020Approximation}.

\begin{lemma}\label{lem:v-apriori}
For $\tilde h\in H^{\frac32}(\partial\Omega)$ and $F\in L^2(\Omega)$, the solution $v_\sigma$ {\rm(}of the penalized problem \eqref{app prob}{\rm)} satisfies
\begin{equation*}
    \|v_\sigma\|_{H^2(\Omega)}\leq C(\|\tilde h\|_{H^\frac32(\partial\Omega)}+\|F\|_{L^2(\Omega)}).
\end{equation*}
\end{lemma}
\begin{proof}
The given conditions on $\tilde h$ and $F$ and the standard elliptic regularity theory imply
$$\|v^*\|_{H^2(\Omega)}\leq C(\|\tilde h\|_{H^\frac32(\partial\Omega)}+\|F\|_{L^2(\Omega)}).$$
Meanwhile, $z=v_\sigma-v^*$ satisfies problem \eqref{eqn:z} with $£\epsilon=\sigma^{-1}$ and $g=-\epsilon \kappa \frac{\partial v^*}{\partial n}$. Then applying Lemma \ref{lemma:Costabel} and the trace theorem yields
\begin{equation*}
    \|v^*-v_\sigma\|_{H^2(\Omega)} \leq C\epsilon^{-1}\|\epsilon\kappa \tfrac{\partial v^*}{\partial n}\|_{H^\frac12(\partial\Omega)} \leq C\|v^*\|_{H^2(\Omega)}.
\end{equation*}
Then the triangle inequality shows the desired assertion.
\end{proof}

The next result gives an error bound on the NN approximation $\widehat{v_{\sigma}}$, by suitably adapting the estimates given in \cite[Theorems 4.1 and 5.13]{JiaoLai:2021elliptic}, and the \textit{a priori} regularity on $v_\sigma$ in Lemma \ref{lem:v-apriori}.
The expectation $\mathbb{E}_{\{X_i\}_{i=1}^N,\{Y_j\}_{j=1}^N}[\cdot]$ is taken with respect to the random sampling points $\{X_i\}_{i=1}^N$ and $\{Y_j\}_{j=1}^N$ in the domain $\Omega$ and on the boundary $\partial\Omega$. We recall also the set $\mathcal{N}_\rho(D,N_\theta,R)$ of NN functions (of a fixed architecture) defined in \eqref{eqn:nn-set}.
\begin{theorem}\label{thm:error}
Let $F\in L^\infty(\Omega)$ and $\tilde h\in H^\frac{3}{2}(\partial\Omega)\cap L^\infty(\partial\Omega)$. Let $N$ be the number of training samples on the domain $\Omega$ and the boundary $\partial\Omega$, and $ \rho $ be logistic function $ \frac{1}{1+{\rm e}^{-x}} $ or tanh function $ \frac{{\rm e}^x-{\rm e}^{-x}}{{\rm e}^x+{\rm e}^{-x}} $. Then for any tolerance $ \epsilon>0 $ and $ \mu\in(0,1) $, there exists {a parameterized NN function class
$$ \mathcal{A}=\mathcal{N}_\rho\left(C\ln(d+1),C(d)\epsilon^{-\frac{d}{1-\mu}},C(d)\epsilon^{-\frac{9d+8}{2-2\mu}}\right)
$$ such that
with the number of samples
$N=O(\epsilon^{-C\frac{d\ln(d+1)}{1-\mu}})$ both in the domain $\Omega$ and on the boundary $\partial\Omega$, the NN approximation $\widehat{v_\sigma}\in\mathcal{A}$ satisfies}
\begin{equation}\label{eqn:err-est0}	 \mathbb{E}_{\{X\}_{i=1}^N,\{Y_j\}_{j=1}^N}\left[\|\nabla(\widehat{v_\sigma}-v_\sigma)\|_{L^2(\Omega)}^2+\sigma\|\widehat{v_\sigma}-v_\sigma\||_{L^2(\partial\Omega)}^2\right]\leq C(\Omega,\kappa,F,\tilde h)\sigma\epsilon^2.
\end{equation}
Thus, for the neural network approximation $\widehat{u_\sigma}=\kappa^{-1}f\Phi+\widehat{v_\sigma}$ by SSDRM, there holds
\begin{equation}\label{err}
\mathbb{E}_{\{X_i\}_{i=1}^N,\{Y_j\}_{j=1}^N}\left[\|\widehat{u_\sigma}-u^*\|_{H^1(\Omega)}^2\right]\leq C(\Omega,\kappa,F,\tilde h)\sigma\epsilon^2+c(F,\tilde h)\sigma^{-2}.
\end{equation}
\end{theorem}
\begin{proof}
The stated result follows from several known estimates from
\cite{JiaoLai:2021elliptic}, and relies on estimating the approximation error and statistical error separately. We sketch the proof for the convenience of readers.
First, we bound the approximation error $\inf_{\widehat{v}\in\mathcal{A}}\|\widehat{v}-v_\sigma\|_{H^1(\Omega)}$ by Lemma \ref{lem:v-apriori}:
\begin{align*}
\inf_{\widehat{v}\in\mathcal{A}} \|\widehat{v}-v_\sigma\|_{H^1(\Omega)}&=\|v_\sigma\|_{H^2(\Omega)}\inf_{\widehat{v}\in\mathcal{A}} \left\|\frac{\widehat{v}}{\|v_\sigma\|_{H^2(\Omega)}}-\frac{v_\sigma}{\|v_\sigma\|_{H^2(\Omega)}}\right\|_{H^1(\Omega)}\\
	&=\|v_\sigma\|_{H^2(\Omega)}\inf_{\widehat{v}\in\mathcal{A}} \left\|\widehat{v}-\frac{v_\sigma}{\|v_\sigma\|_{H^2(\Omega)}}\right\|_{H^1(\Omega)}\\
	&\leq C\big(\|\tilde h\|_{H^\frac32(\partial\Omega)}+\|F\|_{L^2(\Omega)}\big)\inf_{\widehat{v}\in\mathcal{A}} \left\|\widehat{v}-\frac{v_\sigma}{\|v_\sigma\|_{H^2(\Omega)}}\right\|_{H^1(\Omega)}.
\end{align*}
By \cite[Theorem 4.1]{JiaoLai:2021elliptic} (or \cite[Proposition 4.8]{I2020Approximation}), there exists a neural network function $v_\rho\in\mathcal{A}=\mathcal{N}_\rho(C\ln(d+1),C(d)\epsilon^{-\frac{d}{1-\mu}},C(d)\epsilon^{-\frac{9d+8}{2-2\mu}}) $ such that
	$$\left\|v_\rho-\frac{v_\sigma}{\|v_\sigma\|_{H^2(\Omega)}}\right\|_{H^1(\Omega)}\leq\epsilon.$$
Consequently, we get \begin{equation}\label{appro error}\inf_{\widehat{v}\in\mathcal{A}} \|\widehat{v}-v_\sigma\|_{H^1(\Omega)}^2\leq C(\Omega,F,\tilde h)\epsilon^2.
\end{equation}
To bound the statistical error $\sup_{\widehat{v}\in\mathcal{A}}|L_\sigma(\widehat{v})-\widehat{L_\sigma}(\widehat{v})|$, we customarily split it into three terms
\begin{align*}
 \sup_{\widehat{v}\in\mathcal{A}}|L_\sigma(\widehat{v})-\widehat{L_\sigma}(\widehat{v})| \leq&  \sup_{\widehat{v}\in\mathcal{A}}
 \Big||\Omega|\mathbb{E}_{X\sim U(\Omega)}\frac{\kappa(X)|\nabla \widehat{v}(X)|^2}{2}-\frac{|\Omega|}{N}\sum_{i=1}^{N}\frac{\kappa(X_i)|\nabla \widehat{v}(X_i)|^2}{2}\Big|\\
 &+ \sup_{\widehat{v}\in\mathcal{A}} \Big| |\Omega|\mathbb{E}_{X\sim U(\Omega)}\widehat{v}(X)F(X)-\frac{|\Omega|}{N}\sum_{i=1}^{N}\widehat{v}(X_i)F(X_i)\Big|\\
&+\frac{\sigma}{2}\sup_{\widehat{v}\in\mathcal{A}}\Big||\partial\Omega|\mathbb{E}_{Y\sim U(\partial\Omega)}[T\widehat{v}(Y)-\tilde{h}(Y)]^2-\frac{|\partial\Omega|}{N}\sum_{j=1}^{N}[T\widehat{v}(Y_j)-\tilde{h}(Y_j)]^2\Big|.
\end{align*}
Then by applying \cite[Theorem 5.13]{JiaoLai:2021elliptic} (with the following setting of the theorem
statement: the depth $\mathcal{D}=C\ln(d+1)$, the number of nonzero parameters $\boldsymbol{n}_\mathcal{D}=C(d)
\epsilon^{-\frac{d}{1-\mu}}$, and the upper bound $B_\theta=C(d)\epsilon^{-\frac{9d+8}{2-2\mu}}$
on each parameter) to the above three terms, we obtain
\begin{equation}\label{stati error}
		 \mathbb{E}_{\{X_i\}_{i=1}^N,\{Y_j\}_{j=1}^N}\Big[\sup_{\widehat{v}\in\mathcal{A}}|L_\sigma(\widehat{v})-\widehat{L_\sigma}(\widehat{v})|\Big]\leq C(\Omega,\kappa)\sigma\epsilon^2,
\end{equation}
when the number $N$ of sampling points in the domain $\Omega$ and on the boundary $\partial\Omega$ is set to
$N=C(d,\Omega)\epsilon^{-C\frac{d\ln(d+1)}{1-\mu}}$. This estimate was proved using Rademacher complexity and Lipschitz continuity of the NN output with respect to the NN parameter vector $\theta$. Last, substituting the estimates
\eqref{appro error} and \eqref{stati error} into Lemma \ref{lem:decom} yields the desired estimate \eqref{eqn:err-est0}. The other estimate \eqref{err} is direct from the construction of the approximation $\widehat{u_\sigma}$, Poincar\'{e} inequality and Theorem \ref{thm:penalization}.
\end{proof}

Note that generally,  the approximation  $\widehat{u_\sigma}$ and the true solution $u^*$ do not necessarily belong to the space $H^1(\Omega)$, and thus an error estimate in the $H^1(\Omega)$ norm cannot be expected. Nonetheless, the error $\widehat{u_\sigma}-u^*$ does converge in $H^1(\Omega)$, since the singularity splitting technique extracts the leading singularity directly. The penalty parameter $\sigma$ is updated using the path-following strategy, and it is increased geometrically by a factor $\eta>1$ in every outer loop. Thus, as $k$ increases, $\sigma_k$ grows exponentially so the error due to the penalization method decays exponentially fast to zero, cf. Theorem \ref{thm:penalization}. However, this comes at the price of increasing the generalization error, as evidenced by the factor $\sigma$ in the first term of the error bound. This error can be partially offset by $\epsilon$, which is in turn controlled by the number of sample points (in the domain and on the boundary) and the complexity of the neural network (width, depth and weights). This short analysis clearly indicates that there might be an optimal value for the penalty parameter $\sigma$ to balance the trade-off between the different sources of the errors. Of course this analysis has left out one very important point in the analysis: it assumes that a global minimizer $\widehat {v_\sigma}$ to the empirical loss $\widehat{L}_\sigma(v_\theta)$ can be found, which is generally not true, due to the highly complex landscape of the empirical loss $\widehat{L}_\sigma(v_\theta)$.

\section{Conclusions}
In this work, we have investigated the efficient neural network solution of the variable coefficient second-order elliptic problems with singular sources, which include point sources, line sources and their combinations. The presence  of the singular source prevents  a direct application of the standard neural network based solvers, e.g., physics informed neural networks, deep Ritz method and weak adversarial networks. We have proposed  a simple modification of the standard deep Ritz method (which is based on the Ritz variational formulation) by splitting the solution into a regular part and a singular part, where the singular part is expressed analytically using the fundamental solution to the Laplace equation. Extensive numerical experiments clearly confirm its efficiency and accuracy, when compared with existing neural network based approaches. Furthermore, we have discussed relevant theoretical issues.

Although the idea of singularity splitting has only been illustrated on the deep Ritz method, it applies equally well to other neural network based approaches, e.g., physics informed neural networks and weak adversarial networks, and potentially can also be very promising for solving elliptic problems with singular sources, due to the much improved regularity of the regular part of the splitting. These generalizations and related issues, e.g., evolution equations, will be explored in future works.

\section*{Acknowledgements}
The authors are grateful to the two anonymous referees for their many constructive comments which have led to a significant improvement of the quality of the paper.

\bibliographystyle{abbrv}
\bibliography{refer}
	
\end{document}